\documentclass[11pt]{amsart}
\usepackage{amsmath}
\usepackage{amsxtra}
\usepackage{amscd}
\usepackage{amsthm}
\usepackage{amsfonts}
\usepackage{amssymb}
\usepackage{eucal}

\newtheorem{thm}{Theorem}[section]
\newtheorem{conj}[thm]{Conjecture}
\newtheorem{cor}[thm]{Corollary}
\newtheorem{lem}[thm]{Lemma}
\newtheorem{prop}[thm]{Proposition}

\theoremstyle{remark}
\newtheorem{remark}[thm]{Remark}

\theoremstyle{definition}
\newtheorem{definition}[thm]{Definition}

\numberwithin{equation}{section}

\newcommand{\Ref}[1]{{$($\ref{#1}$)$}}
\newcommand{\bean}{\begin{eqnarray}}
\newcommand{\eean}{\end{eqnarray}}
\newcommand{\be}{\begin{displaymath}}
\newcommand{\ee}{\end{displaymath}}
\newcommand{\bea}{\begin{eqnarray*}}   
\newcommand{\eea}{\end{eqnarray*}}

\newcommand{\thmref}[1]{Theorem~\ref{#1}}
\newcommand{\secref}[1]{Section~\ref{#1}}
\newcommand{\lemref}[1]{Lemma~\ref{#1}}
\newcommand{\propref}[1]{Proposition~\ref{#1}}
\newcommand{\corref}[1]{Corollary~\ref{#1}}

\newcommand{\nc}{\newcommand}
\nc{\on}{\operatorname}
\nc{\ch}{\mbox{ch}}
\nc{\Z}{{\mathbb Z}}
\nc{\C}{{\mathbb C}}
\nc{\pone}{{\mathbb C}{\mathbb P}^1}
\nc{\pa}{\partial}
\nc{\F}{{\mathcal F}}
\nc{\arr}{\rightarrow}
\nc{\larr}{\longrightarrow}
\nc{\al}{\alpha}
\nc{\ri}{\rangle}
\nc{\lef}{\langle}
\nc{\W}{{\mathcal W}}
\nc{\la}{\lambda}
\nc{\ep}{\epsilon}
\nc{\su}{\widehat{{\mathfrak sl}}_2}
\nc{\sw}{{\mathfrak s}{\mathfrak l}}
\nc{\g}{{\mathfrak g}}
\nc{\h}{{\mathfrak h}}
\nc{\n}{{\mathfrak n}}
\nc{\N}{\widehat{\n}}
\nc{\G}{\widehat{\g}}
\nc{\De}{\Delta_+}
\nc{\gt}{\widetilde{\g}}
\nc{\Ga}{\Gamma}
\nc{\one}{{\mathbf 1}}
\nc{\z}{{\mathfrak Z}}
\nc{\zz}{{\mathcal Z}}
\nc{\Hh}{{\mathcal H}_\beta}
\nc{\qp}{q^{\frac{k}{2}}}
\nc{\qm}{q^{-\frac{k}{2}}}
\nc{\La}{\Lambda}
\nc{\wt}{\widetilde}
\nc{\wh}{\widehat}
\nc{\qn}{\frac{[m]_q^2}{[2m]_q}}
\nc{\cri}{_{\on{cr}}}
\nc{\kk}{h^\vee}
\nc{\sun}{\widehat{\sw}_N}
\nc{\hh}{{\mathbf H}_{q,t}}
\nc{\HH}{{\mathcal H}_{q,t}}
\nc{\hhh}{{\mathcal H}_{q,1}}
\nc{\ca}{\wt{{\mathcal A}}_{h,k}(\sw_2)}
\nc{\si}{\sigma}
\nc{\gl}{\widehat{{\mathfrak g}{\mathfrak l}}_2}
\nc{\el}{\ell}
\nc{\s}{T}
\nc{\bi}{\bibitem}
\nc{\om}{\omega}
\nc{\WW}{\W_\beta}
\nc{\scr}{{\mathbf S}}
\nc{\ab}{{\mathbf a}}
\nc{\rr}{r}
\nc{\ol}{\overline}
\nc{\con}{qt^{-1} + q^{-1}t}
\nc{\den}{q^{\el-1} t^{-\el+1}+ q^{-\el+1} t^{\el-1}}
\nc{\ds}{\displaystyle}
\nc{\B}{B}
\nc{\A}{A^{(2)}_{2\el}}
\nc{\GG}{{\mathcal G}}
\nc{\UU}{{\mathcal U}}
\nc{\MM}{{\mathcal M}}
\nc{\CC}{{\mathcal C}}
\nc{\GL}{^L\G}
\nc{\gL}{^L\g}
\nc{\dzz}{\frac{dz}{z}}
\nc{\Res}{\on{Res}}
\nc{\rep}{{\mathcal R}ep \;}
\nc{\uqg}{U_q \G}
\nc{\uqgg}{U_q \g}
\nc{\mc}{\mathcal}
\nc{\Cal}{\mathcal}
\nc{\calp}{{\Cal P}}
\nc{\bp}{{\mathbf P}}
\nc{\bq}{{\mathbf Q}}
\nc{\bb}{{\mathfrak b}}
\nc{\uqb}{U_q \bb_-}
\nc{\uqn}{U_q \wt{{\mathfrak n}}}
\nc{\uqh}{U_q \wt{{\mathfrak h}}}
\nc{\uqhh}{U_q \wh{{\mathfrak h}}}
\nc{\uqnn}{U_q \wh{{\mathfrak n}}}
\nc{\ot}{\otimes}
\nc{\R}{{\mc R}}
\nc{\uqbb}{U_q \wt{\g}}
\nc{\yy}{{\mc Y}}
\nc{\uqsl}{U_q \widehat{\sw}_2}
\nc{\ga}{\gamma}
\nc{\Ab}{{\mathbf A}}
\nc{\Yb}{{\mathbf Y}}
\nc{\yb}{{\mathbf y}}
\nc{\uh}{U \wt{{\mathfrak h}}}
\nc{\uhh}{U \wh{{\mathfrak h}}}

\nc{\us}{\underset}
\nc{\opl}{\oplus}
\nc{\yyy}{\wh{\yy}}
\nc{\ovl}{\overline}
\nc{\beq}{\begin{equation}}
\nc{\Fq}{{\mathcal F}}
\nc{\Mq}{{\mathcal M}}
\nc{\Rep}{\on{Rep}}
\nc{\sssec}{\subsubsection}
\nc{\ssec}{\subsection}
\nc{\lan}{\langle}
\nc{\ran}{\rangle}

\nc{\uqhJ}{U_q \widehat{\h}_J^\perp}
\nc{\uqsli}{U_{q_i} \widehat{\sw}_2}
\nc{\len}{\on{length}}

\begin{document}

\title[Combinatorics of the $q$--characters]{Combinatorics of
$q$--characters of finite-dimensional representations of quantum
affine algebras}

\author[Edward Frenkel]{Edward Frenkel}

\author{Evgeny Mukhin}

\address{Department of Mathematics, University of California, Berkeley, CA
94720, USA}

\begin{abstract}
We study finite-dimensional representations of quantum affine algebras
using $q$--characters. We prove the conjectures from \cite{FR:char}
and derive some of their corollaries. In particular, we prove that the
tensor product of fundamental representations is reducible if and only
if at least one of the corresponding normalized $R$--matrices has a
pole.
\end{abstract}

\maketitle

\section*{Introduction}

The intricate structure of the finite-dimensional representations of
quantum affine algebras has been extensively studied from different
points of view, see, e.g.,
\cite{CP,CP3,CP4,CP5,GV,V,KS,AK,FR:char}. While a lot of progress has
been made, many basic questions remained unanswered. In order to
tackle those questions, E. Frenkel and N. Reshetikhin introduced in
\cite{FR:char} a theory of $q$--characters for these
representations. One of the motivations was the theory of deformed
$\W$--algebras developed in \cite{FR:simple}: the representation ring
of a quantum affine algebra should be viewed as a deformed
$\W$--algebra, while the $q$--character homomorphism should be viewed
as its free field realization. The study of $q$--characters in
\cite{FR:char} was based on two main conjectures. One of the goals of
the present paper is to prove these conjectures and to derive some of
their corollaries.

Let us describe our results in more detail. Let $\g$ be a simple Lie
algebra, $\G$ be the corresponding non-twisted affine Kac-Moody
algebra, and $\uqg$ be its quantized universal enveloping algebra
(quantum affine algebra for short). Denote by $I$ the set of vertices
of the Dynkin diagram of $\g$. Let $\on{Rep} \uqg$ be the Grothendieck
ring of $\uqg$. The $q$--character homomorphism is an injective
homomorphism $\chi_q$ from $\on{Rep} \uqg$ to the ring of Laurent
polynomials in infinitely many variables $\yy = \Z[Y_{i,a}^{\pm
1}]_{i\in I; a \in \C^\times}$. This homomorphism should be viewed as
a $q$--analogue of the ordinary character homomorphism.

Indeed, let $G$ be the connected simply-connected algebraic group
corresponding to $\g$, and let $T$ be its maximal torus. We have a
homomorphism $\chi: \on{Rep} G \arr \on{Fun} T$ (where $\on{Fun} T$
stands for the ring of regular functions on $T$), defined by the
formula $(\chi(V))(t) = \on{Tr}_V t$, for all $t \in T$. Upon the
identification of $\on{Rep} G$ with $\on{Rep} \uqgg$ and of $\on{Fun}
T$ with $\Z[y_i^{\pm 1}]_{i\in I}$, where $y_i$ is the function on $T$
corresponding to the fundamental weight $\om_i$, we obtain a
homomorphism $\chi: \on{Rep} \uqgg \arr \Z[y_i^{\pm 1}]_{i\in I}$. One
of the properties of $\chi_q$ is that if we replace each $Y_{i,a}^{\pm
1}$ by $y_i^{\pm 1}$ in $\chi_q(V)$, where $V$ is a $\uqg$--module,
then we obtain $\chi(V|_{\uqgg})$.

The two conjectures from \cite{FR:char} that we prove in this paper
may be viewed as $q$--analogues of the well-known properties of the
ordinary characters. The first of them, \thmref{1}, is the
analogue of the statement that the character of any irreducible
$\uqgg$--module $W$ equals the sum of terms which correspond to the
weights of the form $\la - \sum_{i \in I} n_i \al_i, n_i \in \Z_+$,
where $\la = \sum_{i \in I} l_i \om_i, l_i \in \Z_+$, is the highest
weight of $V$, and $\al_i, i \in I$, are the simple roots. In other
words, we have: $\chi(W) = m_+(1 + \sum_p M_p)$, where $m_+ = \prod_{i
\in I} y_i^{l_i}$, and each $M_p$ is a product of factors $a_j^{-1}, j
\in I$, corresponding to the negative simple roots. \thmref{1} says
that for any irreducible $\uqg$--module $V$, $\chi_q(V) = m_+(1 +
\sum_p M_p)$, where $m_+$ is a monomial in $Y_{i,a}, i \in I, a \in
\C^\times$, with positive powers only (the highest weight monomial),
and each $M_p$ is a product of factors $A_{j,c}^{-1}, j \in I, c \in
\C^\times$, which are the $q$--analogues of the negative simple roots
of $\g$.

The second statement, \thmref{main}, gives an explicit description of
the image of the $q$--character homomorphism $\chi_q$. This is a
generalization of the well-known fact that the image of the ordinary
character homomorphism $\chi$ is equal to the subring of invariants of
$\Z[y_i^{\pm 1}]_{i\in I}$ under the action of the Weyl group $W$ of
$\g$.

Recall that the Weyl group is generated by the simple reflections
$s_i, i \in I$. The subring of invariants of $s_i$ in $\Z[y_i^{\pm
1}]_{i\in I}$ is equal to
$$
K_i = \Z[y_j^{\pm 1}]_{j \neq i} \otimes \Z[y_i + y_i a_i^{-1}],
$$
and hence we obtain a ring isomorphism $\on{Rep} \uqg \simeq \ds
\bigcap_{i\in I} K_i$.

In \thmref{main} (see also \corref{subring}) we establish a
$q$--analogue of this isomorphism. Instead of the simple reflections
we have the screening operators $S_i, i \in I$, introduced in
\cite{FR:char}. We show that $\on{Im} \chi_q$ equals $\ds \bigcap_{i
\in I} \on{Ker} S_i$. Moreover, $\on{Ker} S_i$ is equal to
$$
{\mc K}_i = \Z[Y_{j,a}^{\pm 1}]_{j\neq i; a \in \C^\times} \otimes
\Z[Y_{i,b} + Y_{i,b} A_{i,bq_i}^{-1}]_{b \in \C^\times}.
$$
Thus, we obtain a ring isomorphism $\on{Rep} \uqg \simeq \ds
\bigcap_{i \in I} {\mc K}_i$.

These results allow us to construct in a purely combinatorial way the
$q$--characters of the fundamental representations of $\uqg$, see
\secref{algorithm}.

We derive several corollaries of these results. Here is one of them
(see \thmref{poles} and \propref{irreducible tensor product}). For
each fundamental weight $\omega_i$, there exists a family of
$\uqg$--modules, $V_{\om_i}(a), a \in \C^\times$ (see \secref{fd} for
the precise definition). These are irreducible finite-dimensional
representations of $\uqg$, which have highest weight $\om_i$ if
restricted to $\uqgg$. They are called the fundamental representations
of $\uqg$ (of level $0$). According to a theorem of Chari-Pressley
\cite{CP,CP4} (see \corref{generated} below), any irreducible
representation of $\uqg$ can be realized as a subquotient of a tensor
product of the fundamental representations. The following theorem,
which was conjectured, e.g., in \cite{AK}, describes under what
conditions such a tensor product is reducible.

Denote by $h^\vee$ the dual Coxeter number of $\g$, and by $r^\vee$
the maximal number of edges connecting two vertices of the Dynkin
diagram of $\g$. For the definition of the normalized $R$--matrix, see
\secref{normR}.

\medskip

\noindent{\bf Theorem.} {\em Let $\{ V_k \}_{k=1,\dots,n}$, where $V_k
=V_{\om_{s(k)}}(a_k)$, be a set of fundamental representations of
$\uqg$.}

{\em The tensor product $V_1 \otimes \ldots \otimes V_n$ is reducible
if and only if for some $i,j\in \{ 1,\ldots,n \}$, $i\neq j$, the
normalized $R$--matrix $\ol{R}_{V_i,V_j}(z)$ has a pole at
$z=a_j/a_i$.

In that case $a_j/a_i$ is necessarily equal to $q^k$, where $k$ is an
integer, such that $2\leq |k|\leq r^\vee h^\vee$.}

\medskip

The paper is organized as follows. In Section 1 we recall the main
definitions and results on quantum affine algebras and their
finite-dimensional representations. In Section 2 we give the
definition of the $q$--character homomorphism and list some of its
properties. In Section 3 we develop our main technical tool: the
restriction homomorphisms $\tau_J$. Sections 4 and 5 contain the
proofs of Conjectures 1 and 2 from \cite{FR:char}, respectively. In
Section 6 we use these results to describe the structure of the
$q$--characters of the fundamental representations and to prove the
above Theorem.

The results of this paper can be generalized to the case of the
twisted quantum affine algebras.

In the course of writing this paper we were informed by H. Nakajima
that he obtained an independent proof of Conjecture 1 from
\cite{FR:char} in the $ADE$ case using a geometric approach.

\medskip

\noindent{\bf Acknowledgments.} We thank N.~Reshetikhin for useful
discussions. The research of both authors was supported by a grant
from the Packard Foundation.

\section{Preliminaries on finite-dimensional representations of $\uqg$}

\subsection{Root data}    \label{cartan}

Let $\g$ be a simple Lie algebra of rank $\el$. Let $h^\vee$ be the
dual Coxeter number of $\g$. Let $\langle \cdot,\cdot \rangle$ be the
invariant inner product on $\g$, normalized as in \cite{Kac}, so that
the square of the length of the maximal root equals $2$ with respect
to the induced inner product on the dual space to the Cartan
subalgebra $\h$ of $\g$ (also denoted by $\langle \cdot,\cdot
\rangle$). Denote by $I$ the set $\{ 1,\ldots,\el \}$. Let $\{ \al_i
\}_{i \in I }$ and $\{ \om_i \}_{i \in I}$ be the sets of simple roots
and of fundamental weights of $\g$, respectively. We have:
$$
\langle \al_i,\om_j \rangle = \frac{\langle \al_i,\al_i \rangle}{2}
\delta_{ij}.
$$
Let $r^\vee$ be the maximal number of edges connecting two vertices of
the Dynkin diagram of $\g$. Thus, $r^\vee=1$ for simply-laced $\g$,
$r^\vee=2$ for $B_\el, C_\el, F_4$, and $r^\vee=3$ for $G_2$.

In this paper we will use the rescaled inner product
$$
(\cdot,\cdot) = r^\vee \langle \cdot,\cdot \rangle
$$
on $\h^*$. Set
$$
D = \on{diag}(\rr_1,\ldots,\rr_\el),
$$
where
\begin{equation}    \label{di}
\rr_i = \frac{(\al_i,\al_i)}{2} = r^\vee \frac{\langle \al_i,\al_i
\rangle}{2}.
\end{equation}
The $\rr_i$'s are relatively prime integers. For simply-laced $\g$,
all $r_i$'s are equal to $1$ and $D$ is the identity matrix.

Now let $C = (C_{ij})_{1\leq i,j\leq \el}$
be the {\em Cartan matrix} 
of $\g$, 
$$
C_{ij} = \frac{2(\al_i,\al_j)}{(\al_i,\al_i)}.
$$
Let $B = (B_{ij})_{1\leq i,j\leq \el}$ be the symmetric matrix
$$B = D C,$$ i.e.,
$
B_{ij} = (\al_i,\al_j) = r^\vee \langle \al_i,\al_j \rangle.
$

Let $q \in \C^\times$ be such that $|q|<1$. Set 
$q_i = q^{r_i}$, and $$[n]_q = \frac{q^n - q^{-n}}{q - q^{-1}}.$$
Following \cite{FR:simple,FR:char}, define the $\el \times \el$
matrices $B(q), C(q), D(q)$ by the formulas
\begin{align*}
B_{ij}(q) &= [B_{ij}]_q,\\
C_{ij}(q) &= (q_i + q_i^{-1}) \delta_{ij} + (1-\delta_{ij})[C_{ij}]_q,\\
D_{ij}(q) &= [D_{ij}]_q = \delta_{ij} [r_i]_q.
\end{align*}
We have:
$$
B(q) = D(q) C(q).
$$

Let $\wt{C}(q)$ be the inverse of the Cartan matrix $C(q)$, $C(q)
\wt{C}(q) = \on{Id}$. We will need the following property of matrix
$\wt C(q)$.

\begin{lem}\label{inverse lemma}
All coefficients of the matrix $\wt C(q)$ can be written in the form
\bean\label{inverse formula} \wt C_{ij}(q)=\frac{\wt
C'_{ij}(q)}{d(q)}, \qquad i,j\in I, \eean where $\wt C'_{ij}(q)$,
$d(q)$ are Laurent polynomials in $q$ with non-negative integral
coefficients, symmetric with respect to the substitution $q \arr
q^{-1}$. Moreover, \be \deg \wt{C}'_{ij}(q) < \deg d(q), \qquad
i,j\in I.\ee
\end{lem}

\begin{proof}
We write here the minimal choice of $d(q)$, which we use in
\secref{homtauJ}: \bea A_\el:\; d(q)&=&
q^\el+q^{\el-2}+\dots+q^{-\el}, \\ B_\el:\;
d(q)&=&q^{2\el-1}+q^{2\el-3}+\dots+q^{-2\el-1}, \\ C_\el:\;
d(q)&=&q^{\el+1}+q^{-\el-1}, \\ D_\el:\; d(q) &=&
(q+q^{-1})(q^{\el-1}+q^{-\el+1}), \\ E_6:
\;d(q)&=&(q^2+1+q^{-2})(q^6+q^{-6}), \\ E_7:\;
d(q)&=&(q+q^{-1})(q^9+q^{-9}), \\ E_8:
\;d(q)&=&(q+q^{-1})(q^{15}+q^{-15}), \\ F_4:\; d(q)&=&q^9+q^{-9}, \\
G_2:\; d(q)&=&q^6+q^{-6}.  \eea For Lie algebras of classical series,
the statement of the lemma with the above $d(q)$ follows from the
explicit formulas for the entries $\wt C_{ij}(q)$ of the matrix
$\wt{C}(q)$ given in Appendix C of \cite{FR:simple}. For exceptional
types, the lemma follows from a case by case inspection of the matrix
$\wt C(q)$.
\end{proof}

\subsection{Quantum affine algebras}    \label{qaal}

The quantum affine algebra $\uqg$ in the Drinfeld-Jimbo
realization \cite{Dr1,J} is an associative algebra over $\C$
with generators $x_i^{{}\pm{}}$, $k_i^{{}\pm 1}$ ($i=0,\ldots,\el$),
and relations:
\begin{align*}
k_ik_i^{-1} = k_i^{-1}k_i &=1,\quad \quad k_ik_j =k_jk_i,\\
k_ix_j^{{}\pm{}}k_i^{-1} &= q^{{}\pm B_{ij}}x_j^{{}\pm},\\ [x_i^+ ,
x_j^-] &= \delta_{ij}\frac{k_i - k_i^{-1}}{q_i -q_i^{-1}},\\
\sum_{r=0}^{1-C_{ij}}(-1)^r\left[\begin{array}{cc} 1-C_{ij} \\ r \end{array}
\right]_{q_i}
(x_i^{{}\pm{}})^rx_j^{{}\pm{}}&(x_i^{{}\pm{}})^{1-C_{ij}-r} =0, \ \ \
\ i\ne j.
\end{align*}

Here $(C_{ij})_{0\leq i,j\leq \el}$ denotes the Cartan matrix of $\G$.

The algebra $\uqg$ has a structure of a Hopf algebra with the
comultiplication $\Delta$ and the antipode $S$ given on the generators
by the formulas: \bea \Delta(k_i) &=& k_i \otimes k_i, \\
\Delta(x^+_i) &=& x^+_i \otimes 1 + k_i \otimes x^+_i,\\ \Delta(x^-_i)
&=& x^-_i \otimes k_i^{-1} + 1 \otimes x^-_i, \eea \be S(x_i^+) =
-x_i^+ k_i,\qquad S(x_i^-) = - k_i^{-1} x_i^-,\qquad S(k_i^{\pm 1}) =
k_i^{\mp 1}.  \ee

We define a $\Z$-gradation on $\uqg$ by setting:
$\deg x_0^\pm = \pm 1, \deg x_i^\pm = \deg k_i = 0, i \in
I=\{ 1,\ldots,\el\}$.

Denote the subalgebra of $\uqg$ generated by $k_i^{\pm 1}, x_i^+$
(resp., $k_i^{\pm 1}, x_i^-$), $i=0,\ldots,\el$, by $U_q \bb_+$ (resp.,
$U_q \bb_-$).

The algebra $\uqgg$ is defined as the subalgebra of $\uqg$ with
generators $x_i^{{}\pm{}}$, $k_i^{{}\pm 1}$, where $i\in I$.

We will use Drinfeld's ``new'' realization of $\uqg$, see \cite{Dr},
described by the following theorem.

\begin{thm}[\cite{Dr,KhT2,LSS,Beck}]    \label{defining}
The algebra $\uqg$ has another realization as the algebra with
generators $x_{i,n}^{\pm}$ ($i\in I$, $n\in\Z$),
$k_i^{\pm 1}$ ($i\in I$), $h_{i,n}$ ($i\in I$, $n\in
\Z\backslash 0$) and central elements $c^{\pm 1/2}$, with the
following relations:
\begin{align*}
k_ik_j = k_jk_i,\quad & k_ih_{j,n} =h_{j,n}k_i,\\
k_ix^\pm_{j,n}k_i^{-1} &= q^{\pm B_{ij}}x_{j,n}^{\pm},\\
[h_{i,n} , x_{j,m}^{\pm}] &= \pm \frac{1}{n} [n B_{ij}]_q c^{\mp
{|n|/2}}x_{j,n+m}^{\pm},\\ x_{i,n+1}^{\pm}x_{j,m}^{\pm}
-q^{\pm B_{ij}}x_{j,m}^{\pm}x_{i,n+1}^{\pm} &=q^{\pm
B_{ij}}x_{i,n}^{\pm}x_{j,m+1}^{\pm}
-x_{j,m+1}^{\pm}x_{i,n}^{\pm},\\ [h_{i,n},h_{j,m}]
&=\delta_{n,-m} \frac{1}{n} [n B_{ij}]_q \frac{c^n -
c^{-n}}{q-q^{-1}},\\ [x_{i,n}^+ , x_{j,m}^-]=\delta_{ij} & \frac{
c^{(n-m)/2}\phi_{i,n+m}^+ - c^{-(n-m)/2} \phi_{i,n+m}^-}{q_i -
q_i^{-1}},\\
\sum_{\pi\in\Sigma_s}\sum_{k=0}^s(-1)^k\left[\begin{array}{cc} s \\
k \end{array} \right]_{q_i} x_{i, n_{\pi(1)}}^{\pm}\ldots
x_{i,n_{\pi(k)}}^{\pm} & x_{j,m}^{\pm} x_{i,
n_{\pi(k+1)}}^{\pm}\ldots x_{i,n_{\pi(s)}}^{\pm} =0,\ \  s=1-C_{ij},
\end{align*}
for all sequences of integers $n_1,\ldots,n_s$, and $i\ne j$, where
$\Sigma_s$ is the symmetric group on $s$ letters, and
$\phi_{i,n}^{\pm}$'s are determined by the formula
\begin{equation}    \label{series}
\Phi_i^\pm(u) := \sum_{n=0}^{\infty}\phi_{i,\pm n}^{\pm}u^{\pm n} =
k_i^{\pm 1} \exp\left(\pm(q-q^{-1})\sum_{m=1}^{\infty}h_{i,\pm m}
u^{\pm m}\right).
\end{equation}
\end{thm}

For any $a \in \C^\times$, there is a Hopf algebra automorphism
$\tau_a$ of $\uqg$ defined on the generators by the following
formulas:
\bean\label{tau}
\tau_a(x_{i,n}^{\pm})=a^nx_{i,n}^{\pm}, \quad \quad
\tau_a(\phi_{i,n}^{\pm})=a^n\phi_{i,n}^{\pm},
\eean
$$
\tau_a(c^{1/2})=c^{1/2}, \quad \quad \tau_a(k_i)=k_i,
$$
for all $i \in I, n \in \Z$. Given a $\uqg$--module $V$ and $a \in
C^\times$, we denote by $V(a)$ the pull-back of $V$ under $\tau_a$.

Define new variables $\wt{k}_i^{\pm 1}, i\in I$, such that
\begin{equation}    \label{wtk}
k_j = \prod_{i\in I} \wt{k}_i^{C_{ij}}, \quad \quad \wt{k}_i \wt{k}_j
= \wt{k}_j \wt{k}_i.
\end{equation}
Thus, while $k_i$ corresponds to the simple root $\al_i$, $\wt{k}_i$
corresponds to the fundamental weight $\omega_i$. We extend the
algebra $\uqg$ by replacing the generators $k_i^{\pm 1}, i \in I$ with
$\wt{k}_i^{\pm 1}, i\in I$. From now on $\uqg$ will stand for the
extended algebra. 

Let $q^{2\rho} = \wt{k}_1^2 \ldots \wt{k}_\el^2$. The square of
the antipode acts as follows (see \cite{Dr2}):
\bean\label{ss}
S^2(x)=\tau_{q^{-2r^\vee h^\vee}}(q^{-2\rho}x q^{2\rho}), \qquad
\forall x\in\uqg.
\eean

Let $w_0$ be the longest element of the Weyl group of ${\mathfrak
g}$. Let $i\to\bar i$ be the bijection $I\to I$, such that
$w_0(\al_i)=-\al_{\bar i}$. Define the algebra automorphism
$w_0:\uqg\to\uqg$ by \bean\label{w0} w_0(\wt{k}_{i}) = \wt{k}_{\bar
i},\qquad w_0(h_{i,n})=h_{\bar i,n},\qquad
w_0(x^{\pm}_{i,n})=x^\pm_{\bar i,n}.  \eean We have: $w_0^2 =
\on{Id}$. Actually, $w_0$ is a Hopf algebra automorphism, but we will
not use this fact.

\subsection{Finite-dimensional representations of $\uqg$}
\label{fd}

In this section we recall some of the results of Chari and Pressley
\cite{CP,CP3,CP4,CP5} on the structure of finite-dimensional
representations of $\uqg$.

Let $P$ be the weight lattice of $\g$. It is equipped with the
standard {\em partial order}: the weight $\la$ is higher than the
weight $\mu$ if $\la-\mu$ can be written as a combination of the
simple roots with positive integral coefficients.

A vector $w$ in a $\uqgg$--module $W$ is called a vector of weight
$\lambda \in P$, if
\begin{equation}    \label{weight}
k_i \cdot w = q^{(\lambda,\al_i)} w, \quad \quad i \in I.
\end{equation}
A representation $W$ of $\uqgg$ is said to be of type 1 if it is the
direct sum of its weight spaces $W = \ds \oplus_{\la \in P} W_\la$,
where $W_{\lambda} = \{w\in W | k_i \cdot w = q^{(\lambda,\al_i)}w\}$.
If $W_\lambda\ne 0$, then $\lambda$ is called a weight of $W$.

A representation $V$ of $\uqg$ is called of type 1 if $c^{1/2}$ acts
as the identity on $V$, and if $V$ is of type 1 as a representation of
$\uqgg$. According to \cite{CP}, every finite-dimensional irreducible
representation of $\uqg$ can be obtained from a type 1 representation
by twisting with an automorphism of $\uqg$. Because of that, we will
only consider type 1 representations in this paper.

A vector $v\in V$ is called a {\em highest weight vector} if
\begin{equation}    \label{hwv}
x_{i,n}^+ \cdot v=0,\quad \quad \phi_{i,n}^\pm \cdot v =
\psi_{i,n}^\pm v,\quad \quad c^{1/2} v =v, \quad \quad \forall i \in
I, n \in \Z,
\end{equation}
for some complex numbers $\psi_{i,n}^{\pm}$. A type 1 representation
$V$ is a {\em highest weight representation} if $V=\uqg \cdot v$, for
some highest weight vector $v$. In that case the set of generating functions
$$
\Psi_i^\pm(u) = \sum_{n=0}^\infty \psi_{i,\pm n}^\pm u^{\pm n}, \quad
\quad i \in I,
$$
is called the {\em highest weight} of $V$.

\medskip

\noindent{\em Warning.} The above notions of highest weight vector and
highest weight representation are different from standard. Sometimes
they are called pseudo-highest weight vector and pseudo-highest weight
representation.

\medskip

Let $\calp$ be the set of all $I$--tuples $(P_i)_{i\in I}$ of
polynomials $P_i\in\C[u]$, with constant term 1.

\begin{thm}[\cite{CP,CP4}]    \label{ChP2}
\hfill

(1) Every finite-dimensional irreducible representation of $\uqg$ of
type 1 is a highest weight representation.

(2) Let $V$ be a finite-dimensional irreducible representation of
$\uqg$ of type 1 and highest weight $(\Psi_{i}^{\pm}(u))_{i\in
I}$. Then, there exists $\bp=(P_i)_{i\in I}\in\calp$ such that
\begin{equation}\label{hw}
\Psi_i^\pm(u) = q_i^{\deg(P_i)}\frac{P_i(uq_i^{-1})}{P_i(uq_i)},
\end{equation}
as an element of $\C[[u^{\pm 1}]]$.

Assigning to $V$ the $I$--tuple $\bp \in \calp$ defines a bijection
between $\calp$ and the set of isomorphism classes of
finite-dimensional irreducible representations of $\uqg$ of type
1. The irreducible representation associated to $\bp$ will be denoted
by $V(\bp)$.

(3) The highest weight of $V(\bp)$ considered as a $\uqgg$--module is
$\la = \sum_{i\in I} \deg P_i \cdot \omega_i$, the lowest weight of
$V(\bp)$ is $\ol{\la} = - \sum_{i\in I} \deg P_i \cdot
\omega_{\bar{i}}$, and each of them has multiplicity $1$.

(4) If $\bp=(P_i)_{i\in I}\in\calp$, $a\in\C^\times$, and if
$\tau_a^*(V(\bp))$ denotes the pull-back of $V(\bp)$ by the automorphism
$\tau_a$, we have $\tau_a^*(V(\bp))\cong V(\bp^a)$
as representations of $\uqg$, where $\bp^a=(P_i^a)_{i\in I}$ and
$P_i^a(u)=P_i(ua)$.

(5) For $\bp$, $\bq\in\calp$ denote by $\bp\ot\bq \in \calp$ the
$I$--tuple $(P_iQ_i)_{i\in I}$. Then $V(\bp\ot\bq)$ is isomorphic to a
quotient of the subrepresentation of $V(\bp)\ot V(\bq)$ generated by
the tensor product of the highest weight vectors.
\end{thm}

An analogous classification result for Yangians has been
obtained earlier by Drinfeld \cite{Dr}. Because of that, the
polynomials $P_i(u)$ are called Drinfeld polynomials.

Note that in our notation the polynomials $P_i(u)$
correspond to the polynomials $P_i(uq_i^{-1})$ in the notation of
\cite{CP,CP4}.

For each $i\in I$ and $a \in \C^\times$, define the irreducible
representation $V_{\omega_i}(a)$ as $V(\bp^{(i)}_a)$, where
$\bp^{(i)}_a$ is the $I$--tuple of polynomials, such that
$P_i(u)=1-ua$ and $P_j(u)=1, \forall j \neq i$. We call
$V_{\omega_i}(a)$ the $i$th {\em fundamental representation} of
$\uqg$. Note that in general $V_{\omega_i}(a)$ is reducible as a
$\uqgg$--module.

\thmref{ChP2} implies the following

\begin{cor}[\cite{CP4}]    \label{generated}
Any irreducible finite-dimensional representation $V$ of $\uqg$ occurs
as a quotient of the submodule of the tensor product
$V_{\omega_{i_1}}(a_1) \otimes \ldots \otimes V_{\omega_{i_n}}(a_n)$,
generated by the tensor product of the highest weight vectors. The
parameters $(\omega_{i_k},a_k)$, $k=1,\dots,n$, are uniquely
determined by $V$ up to permutation.
\end{cor}

\section{Definition and first properties of $q$--characters}

\subsection{Definition of $q$--characters}    \label{defin}

Let us recall the definition of the $q$--characters of
finite-dimensional representations of $\uqg$ from \cite{FR:char}.

The completed tensor product $\uqg \;\widehat{\otimes}\; \uqg$ contains a
special element $\R$ called the universal $R$--matrix (at level
$0$). It actually lies in $U_q \bb_+ \;\widehat{\otimes}\; U_q \bb_-$ and
satisfies the following identities:
\begin{align*}
\Delta'(x) &= \R \Delta(x) \R^{-1}, \quad \quad \forall x
\in \uqg,\\
(\Delta \otimes \on{id}) \R &= \R^{13} \R^{23}, \quad
\quad (\on{id} \otimes \Delta) \R = \R^{13} \R^{12}.
\end{align*}
For more details, see \cite{Dr2,EFK}.

Now let $(V,\pi_V)$ be a finite-dimensional representation of
$\uqg$. Define the transfer-matrix corresponding to $V$ by
\begin{equation}    \label{tv}
t_V = t_V(z) = \on{Tr}_V \; (\pi_{V(z)} \otimes
\on{id})(\R).
\end{equation}
Thus we obtain a map $\nu_q: \on{Rep} \uqg \arr \uqb[[z]]$, sending
$V$ to $t_V(z)$.

\begin{remark} Note that in \cite{FR:char} there was an extra factor
$q^{2\rho}$ in formula \eqref{tv}. This factor is inessential for the
purposes of this paper, and therefore can be dropped.\qed
\end{remark}

Denote by $\uqbb$ the subalgebra of $\uqg$ generated by $x^\pm_{i,n},
\wt{k}_i, h_{i,r}, n \leq 0, r<0, i\in I$. It follows from the proof of
\thmref{defining} that $U_q\bb_- \subset \uqbb$. As a vector space,
$\uqbb$ can be decomposed as follows: $\uqbb = \uqn_- \otimes \uqh
\otimes \uqn_+$, where $\uqn_\pm$ (resp., $\uqh$) is generated by
$x^\pm_{i,n}, i \in I, n \leq 0$ (resp., $\wt{k}_i, h_{i,n}, i \in I,
n<0$). Hence
\begin{equation*}
\uqbb = \uqh \oplus \left( \uqbb \cdot (\uqn_+)_0 + (\uqn_-)_0
\cdot \uqbb \right),
\end{equation*}
where $(\uqn_\pm)_0$ stands for the augmentation ideal of
$\uqn_\pm$. Denote by ${\mathbf h}_q$ the projection $\uqbb \arr \uqh$
along the last two summands (this is an analogue of the Harish-Chandra
homomorphism). We denote by the same letter its restriction to $\uqb$.

Now we define the map $\chi_q: \on{Rep} \uqg \arr \uqh[[z]]$ as the
composition of $\nu_q: \on{Rep} \uqg \arr \uqb[[z]]$ and ${\mathbf
h}_q[[z]]: \uqb[[z]] \arr \uqh[[z]]$.

To describe the image of $\chi_q$ we need to introduce some more
notation.

Let
\begin{equation}    \label{connection}
\wt{h}_{i,m} = \sum_{j\in I} \wt{C}_{ji}(q^m) h_{j,m},
\end{equation}
where $\wt{C}(q)$ is the inverse matrix to $C(q)$ defined in
\secref{cartan}. Set
\begin{equation}    \label{Yia}
Y_{i,a} = \wt{k}_i^{-1} \exp \left( - (q-q^{-1}) \sum_{n>0}
\wt{h}_{i,-n} z^n a^n \right), \quad \quad a \in \C^\times.
\end{equation}
We assign to $Y_{i,a}^{\pm 1}$ the weight $\pm \omega_i$.

We have the ordinary character homomorphism $\chi: \on{Rep} \uqgg \arr
\Z[y_i^{\pm 1}]_{i\in I}$: if $V = \oplus_\mu V_\mu$ is the weight
decomposition of $V$, then $\chi(V) = \sum_\mu \dim V_\mu \cdot
y^\mu$, where for $\mu = \sum_{i\in I} m_i \om_i$ we set $y^\mu =
\prod_{i\in I} y_i^{m_i}$. Define the homomorphism
$$\beta: \Z[Y_{i,a}^{\pm 1}]_{i\in I; a \in \C^\times} \arr
\Z[y_i^{\pm 1}]_{i\in I}$$ sending $Y_{i,a}^{\pm 1}$ to $y_i^{\pm 1}$,
and denote by
$$
\on{res}: \on{Rep} \uqg \arr \on{Rep} \uqgg
$$
the restriction homomorphism.

Given a polynomial ring $\Z[x_\al^{\pm 1}]_{\al \in A}$, we denote by
$\Z_+[x_\al^{\pm 1}]_{\al \in A}$ its subset consisting of all linear
combinations of monomials in $x_\al^{\pm 1}$ with positive integral
coefficients.

\begin{thm}[\cite{FR:char}]    \label{mainthm}
\hfill

(1) $\chi_q$ is an injective homomorphism from $\on{Rep} \uqg$ to
$\Z[Y_{i,a}^{\pm 1}]_{i\in I; a \in \C^\times} \subset \uqh[[z]]$.

(2) For any finite-dimensional representation $V$ of $\uqg$,
$\chi_q(V) \in \Z_+[Y_{i,a}^{\pm 1}]_{i\in I; a \in \C^\times}$.

(3) The diagram
$$\begin{CD} \on{Rep} \uqg @>{\chi_q}>> \Z[Y_{i,a}^{\pm 1}]_{i\in I; a
\in \C^\times}\\ @VV{\on{res}}V @VV{\beta}V\\ \on{Rep} \uqgg
@>{\chi}>> \Z[y_i^{\pm 1}]_{i\in I}
\end{CD}$$
is commutative.

(4) $\on{Rep} \uqg$ is a commutative ring that is isomorphic to
$\Z[t_{i,a}]_{i\in I; a \in \C^\times}$, where $t_{i,a}$ is the class
of $V_{\omega_i}(a)$.
\end{thm}

The homomorphism
$$
\chi_q: \on{Rep} \uqg \arr \Z[Y_{i,a}^{\pm 1}]_{i\in I; a \in
\C^\times}
$$
is called the $q$--{\em character homomorphism}. For a
finite-dimensional representation $V$ of $\uqg$, $\chi_q(V)$ is called
the $q$--{\em character} of $V$.

\subsection{Spectra of $\Phi^\pm(u)$}

According to \thmref{mainthm}(1), the $q$--character of any
finite-dimensional representation $V$ of $\uqg$ is a linear
combination of monomials in $Y_{i,a}^{\pm 1}$ with positive integral
coefficients. The proof of \thmref{mainthm} from \cite{FR:char} allows
us to relate the monomials appearing in $\chi_q(V)$ to the spectra of
the operators $\Phi_i^\pm(u)$ on $V$ as follows.

It follows from the defining relations that the operators
$\phi_{i,n}^\pm$ commute with each other. Hence we can decompose any
representation $V$ of $\uqg$ into a direct sum $V = \oplus
V_{(\ga^\pm_{i,n})}$ of generalized eigenspaces
$$
V_{(\ga^\pm_{i,n})} = \{ x \in V |\; {\rm there\;\; exists}\;\; p,\;
{\rm such \;\; that}\;\; (\phi^\pm_{i,n} - \ga^\pm_{i,n})^p \cdot x =
0, \forall i \in I,n \in \Z \}.
$$
Since $\phi_0^\pm = k_i^{\pm 1}$, all vectors in
$V_{(\ga^\pm_{i,n})}$ have the same weight (see formula \eqref{weight}
for the definition of weight). Therefore the decomposition of $V$ into
a direct sum of subspaces $V_{(\ga^\pm_{i,n})}$ is a refinement of its
weight decomposition.

Given a collection $(\ga^\pm_{i,n})$ of generalized eigenvalues, we
form the generating functions
$$
\Ga^\pm_i(u) = \sum_{n\geq 0} \ga^\pm_{i,\pm n} u^{\pm n}.
$$
We will refer to each collection $\{ \Ga^\pm_i(u) \}_{i\in I}$
occurring on a given representation $V$ as the {\em common
(generalized) eigenvalues} of $\Phi^\pm_i(u), i\in I$, on $V$, and to
$\dim V_{(\ga^\pm_{i,n})}$ as the {\em multiplicity} of this
eigenvalue.

Let ${\mathfrak B}_V$ be a Jordan basis of $\phi_{i,n}^\pm, i \in I, n
\in \Z$. Consider the module $V(z) = \tau_z^*(V)$, see formula
\Ref{tau}.  Then $V(z)=V$ as a vector space. Moreover, the
decomposition in the direct sum of generalized eigenspaces of
operators $\phi_{i,n}^\pm$ does not depend on $z$, because the action
of $\phi_{i,n}^\pm$ on $V$ and on $V(z)$ differs only by scalar
factors $z^n$. In particular, ${\mathfrak B}_V$ is also a Jordan basis
for $\phi_{i,n}^\pm$ acting on $V(z)$ for all $z\in\C^\times$. If
$v\in{\mathfrak B}_V$ is a generalized eigenvector with common
eigenvalues $\{ \Ga^\pm_i(u) \}_{i\in I}$, then the corresponding
common eigenvalues on $v$ in $V(z)$ are $\{\Ga^\pm_i(zu) \}_{i\in I}$

The following result is a generalization of \thmref{ChP2}.

\begin{prop}[\cite{FR:char}]    \label{spectra}
The eigenvalues $\Ga^\pm_i(u)$ of $\Phi_i^\pm(u)$ on any
finite-dimensional representation of $\uqg$ have the form:
\begin{equation}    \label{form}
\Ga_i^\pm(u) = q_i^{\deg Q_i - \deg R_i} \frac{Q_i(uq_i^{-1})
R_i(uq_i)}{Q_i(uq_i) R_i(uq_i^{-1})},
\end{equation}
as elements of $\C[[u^{\pm 1}]]$, where $Q_i(u), R_i(u)$ are
polynomials in $u$ with constant term $1$.
\end{prop}

Now we can relate the monomials appearing in $\chi_q(V)$ to the common
eigenvalues of $\Phi_i^\pm(u)$ on $V$.

\begin{prop}    \label{common spectra}
Let $V$ be a finite-dimensional $\uqg$--module. There is a one-to-one
correspondence between the monomials occurring in $\chi_q(V)$ and the
common eigenvalues of $\Phi_i^\pm(u), i\in I$, on $V$. Namely, the
monomial
\begin{equation}    \label{monomial}
\prod_{i\in I} \left( \prod_{r=1}^{k_i} Y_{i,a_{ir}} \prod_{s=1}^{l_i}
Y_{i,b_{is}}^{-1} \right)
\end{equation}
corresponds to the common eigenvalues \eqref{form}, where
\begin{equation}    \label{qiri}
Q_i(z) = \prod_{r=1}^{k_i} (1-za_{ir}), \quad \quad R_i(z) =
\prod_{s=1}^{l_i} (1-zb_{is}), \quad \quad i \in I.
\end{equation}
The weight of each monomial equals the weight of the corresponding
generalized eigen\-space. Moreover, the coefficient of each monomial
in $\chi_q(V)$ equals the multiplicity of the corresponding common
eigenvalue.
\end{prop}

\begin{proof}
Denote by $\uqnn_\pm$ the subalgebra of $\uqg$ generated by
$x^\pm_{i,n}, i\in I, n\in \Z$. Let $\wt{B}(q)$ be the inverse matrix
to $B(q)$ from \secref{cartan}. The following formula for the universal
$R$--matrix has been proved in \cite{KhT2,LSS,Da}:
\begin{equation}    \label{ktr}
\R  = \R^+ \R^0 \R^- T,
\end{equation}
where
\begin{equation}
\R^0 = \exp \left( -\sum_{n>0}\sum_{i \in
I}\frac{n(q-q^{-1})^2}{q_i^n-q_i^{-n}} h_{i,n} \otimes \wt{h}_{i,-n}
z^n \right)
\end{equation}
(here we use the notation \eqref{connection}), $\R^\pm \in \uqnn_\pm
\otimes \uqn_\mp$, and $T$ acts as follows: if $x,y$ satisfy $k_i
\cdot x = q^{(\la,\al_i)} x, k_i \cdot y = q^{(\mu,\al_i)} y$, then
\begin{equation}    \label{t}
T \cdot x \otimes y = q^{-(\la,\mu)} x \otimes y.
\end{equation}

By definition, $\chi_q(V)$ is obtained by taking the trace of
$(\pi_{V(z)} \otimes \on{id})(\R)$ over $V$ and then
projecting it on $\uqh[[z]]$ using the projection operator ${\mathbf
h}_q$. This projection eliminates the factor $\R^-$, and then taking
the trace eliminates $\R^+$ (recall that $\uqn_+$ acts nilpotently on
$V$). Hence we obtain:
\begin{equation}\label{explicit character}
\chi_q(V) = \on{Tr}_V\left[ \on{exp} \left(
- \sum_{n>0}\sum_{i \in I}\frac{n(q-q^{-1})^2}{q_i^n-q_i^{-n}}
\pi_V(h_{i,n})\otimes \wt{h}_{i,-n} z^n \right)(\pi_V\otimes 1)T\right],
\end{equation}

The trace can be written as the sum of terms $m_v$ corresponding to
the (generalized) eigenvalues of $h_{i,n}$ on the vectors $v$ of the
Jordan basis ${\mathfrak B}_V$ of $V$ for the operators
$\phi_{i,n}^\pm$ (and hence for $h_{i,n}$).

The eigenvalues of $\Phi_i^\pm(u)$ on each vector $v \in {\mathfrak
B}_V$ are given by formula \eqref{form}. Suppose that $Q_i(u)$ and
$R_i(u)$ are given by formula \eqref{qiri}. Then the eigenvalue of
$h_{i,n}$ on $v$ equals
\begin{equation}    \label{h eig}
\frac{q_i^n-q_i^{-n}}{n(q-q^{-1})} \left( \sum_{r=1}^{k_i}
(a_{ir})^n - \sum_{s=1}^{l_i} (b_{is})^n \right), \quad \quad n>0.
\end{equation}

Substituting into formula \eqref{explicit character} and recalling the
definition \eqref{Yia} of $Y_{i,a}$ we obtain that the corresponding
term $m_v$ in $\chi_q(V)$ is the monomial \eqref{monomial}.
\end{proof}

Let $V = V({\bp})$, where
\begin{equation}    \label{formbp}
P_i(u) = \prod_{k=1}^{n_i} (1-ua^{(i)}_k), \quad \quad i \in I.
\end{equation}
Then by \thmref{ChP2}(3), the module $V$ has
highest weight $\lambda = \sum_{i\in I} \deg P_i \cdot \omega_i$,
which has multiplicity $1$. \propref{common spectra} implies that
$\chi_q(V)$ contains a unique monomial of weight $\lambda$. This
monomial equals
\begin{equation}    \label{mp}
\prod_{i\in I} \prod_{k=1}^{n_i} Y_{i,a^{(i)}_k}.
\end{equation}
We call it the {\em highest weight monomial} of $V$. All other
monomials in $\chi_q(V)$ have lower weight than $\lambda$.

A monomial in $\Z[Y_{i,a}^{\pm 1}]_{i\in I,a\in\C^\times}$ is called
{\em dominant} if it does not contain factors $Y_{i,a}^{-1}$ (i.e., if it
is a product of $Y_{i,a}$'s in positive powers only). The highest
weight monomial is dominant, but in general the highest weight
monomial is not the only dominant monomial occurring in
$\chi_q(V)$. Nevertheless, we prove below in \corref{no dominant
monomials} that the only dominant monomial contained in the
$q$--character of a fundamental representation $V_{\omega_i}(a)$ is
its highest weight monomial $Y_{i,a}$.

Note that a dominant monomial has dominant weight but not all
monomials of dominant weight are dominant.

Similarly, a monomial in $\Z[Y_{i,a}^{\pm 1}]_{i\in I,a\in\C^\times}$
is called {\em antidominant} if it does not contain factors $Y_{i,a}$
(i.e., if it is a product of $Y^{-1}_{i,a}$'s in negative powers
only). The roles of dominant and antidominant monomials are similar,
see, e.g., Remark \ref{dominant-antidominant}.  By Corollary \ref{lowest in
arbitrary}, the lowest weight monomial is antidominant.

\begin{remark}
The statement analogous to \propref{spectra} in the case of the
Yangians has been proved by Knight \cite{Kn}. Using this statement, he
introduced the notion of character of a representation of Yangian.\qed
\end{remark}

\subsection{Connection with the entries of the $R$--matrix}
\label{normR}

We already described the $q$--character of $\uqg$ module $V$ in terms
of universal $R$-matrix and in terms of generalized eigenvalues of
operators $\phi_{i,n}^\pm$. It allows us to describe the
$q$--character of $V$ in terms of diagonal entries of $R$-matrices
acting on the tensor products $V\otimes V_{\om_i}(a)$ with fundamental
representations. We will use this description in Section
\ref{fundamental representations}.

Define
\begin{equation}    \label{Aia}
A_{i,a} = k_i^{-1} \exp \left( - (q-q^{-1}) \sum_{n>0} h_{i,-n}
z^n a^n \right), \quad \quad a \in \C^\times.
\end{equation}
Using formula \eqref{connection}, we can express $A_{i,a}$ in terms of
$Y_{j,b}$'s:
\begin{equation}    \label{express}
A_{i,a} = Y_{i,aq_i} Y_{i,aq_i^{-1}} \prod_{C_{ji}=-1} Y_{j,a}^{-1}
\prod_{C_{ji}=-2} Y_{j,aq}^{-1} Y_{j,aq^{-1}}^{-1} \prod_{C_{ji}=-3}
Y_{j,aq^2}^{-1} Y_{j,a}^{-1} Y_{j,aq^{-2}}^{-1}.
\end{equation}
Thus, $A_{i,a} \in \Z[Y_{j,b}^{\pm 1}]_{j\in I;b\in\C^\times}$, and
the weight of $A_{i,a}$ equals $\al_i$.

Let $V$ and $W$ be irreducible finite-dimensional representations of
$\uqg$ with highest weight vectors $v$ and $w$. Let
$\ol{R}_{VW}(z)\in\on{End} (V\otimes W)$ be the normalized $R$-matrix,
\be \ol{R}_{VW}(z)=f^{-1}_{VW}(z)(\pi_{V(z)}\otimes \pi_W)({\mc R}),
\ee where $f_{VW}(z)$ is the scalar function, such that
\begin{equation}    \label{sends}
\ol{R}_{VW}(z)(v\otimes w) = w\otimes v.
\end{equation}

In what follows we always consider the normalized $R$-matrix
$\ol{R}_{VW}(z)$ written in the basis ${\mathfrak B}_V\otimes {\mathfrak
B}_W$.

Recall the definition of the fundamental representation $V_{\om_i}(a)$
from \secref{fd}. Denote its highest weight vector by $v_{\om_i}$.

\begin{lem}\label{diagonal lemma}
Let $v\in{\mathfrak B}_V$ and suppose that the corresponding monomial
$m_v$ in $\chi_q(V)$ is given by \bean \label{form1} m_v=m_+ M \prod_k
A^{-1}_{i,a_k}, \eean where $M$ is a product of factors
$A_{j,b}^{-1}$, $b\in\C^\times$, $j\in I$, $j\neq i$. Then the
diagonal entry of the normalized $R$-matrix
$\ol{R}_{V,V_{\om_i}(b)}(z)$ corresponding to the vector $v\otimes
v_{\om_i}$ is \bean\label{diagonal}
\left(\ol{R}_{V,V_{\om_i}(b)}(z)\right)_{v\otimes v_{\om_i},v\otimes
v_{\om_i}}= \prod_k q_i \frac{1-a_kzb^{-1}
q_i^{-1}}{1-a_kzb^{-1}q_i}.  \eean
\end{lem}
\begin{proof}
Recall formula \eqref{ktr} for $\R$. We have: $\R^-(v\otimes
v_{\om_i}) = 0$; $v\otimes v_{\om_i}$ is a generalized eigenvector of
$\R^0$; and $\R^+(v\otimes v_{\om_i})$ is a linear combination of
tensor products $x \otimes y \in {\mathfrak B}_V \otimes {\mathfrak
B}_{V_{\om_i}(b)}$, where $y$ has a lower weight than
$v_{\om_i}$. Therefore the diagonal matrix element of $\R$ on $v
\otimes v_{\om_i} \in V(z) \otimes V_{\om_i}(b)$ equals the
generalized eigenvalue of $(\pi_{V(z)} \otimes
\pi_{V_{\om_i}(b)})(\R^0)$ on $v\otimes v_{\om_i}$.

On the other hand, as explained in the proof of \propref{common
spectra}, the monomial $m_v$ is equal to the diagonal matrix element
of $(\pi_{V(z)} \otimes 1)(\R^0)$ corresponding to $v$. Therefore the
diagonal matrix element of $\R$ corresponding to $v \otimes v_{\om_i}$
equals the eigenvalue of $m_v$ (considered as an element of
$\uqh[[z]]$) on $v_{\om_i}$.

In particular, if $v$ is the highest weight vector, then the
corresponding monomial $m_v$ is the highest weight monomial
$m_+$. Therefore we find that the diagonal matrix element of the
non-normalized $R$--matrix corresponding to $v \otimes v_{\om_i}$
equals the eigenvalue of $m_+$ on $v_{\om_i}$. By formula
\eqref{sends} the diagonal matrix element of the normalized
$R$--matrix equals $1$. Therefore the eigenvalue of $m_+$ on
$v_{\om_i}$ equals the scalar function
$f_{V,V_{\om_i}(b)}(z)$. Therefore we obtain that the diagonal matrix
element of the normalized $R$-matrix $\ol{R}_{V,V_{\om_i}(b)}(z)$
corresponding to the vector $v \otimes v_{\om_i}$ is equal to the
eigenvalue of $m_v m_+^{-1}$ on $v_{\om_i}$. According to formula
\eqref{Aia}, $A_{i,a} = \Phi_i^-(za)$. Therefore, if $m_v$ is given by
formula \eqref{form1}, we obtain from formula \eqref{hw} that this
matrix element is given by formula \eqref{diagonal}.
\end{proof}

Note that by Theorem \ref{1} below every monomial occurring in the
$q$--character of an irreducible representation $V$ can be written in
the form \eqref{form1}.

\section{The homomorphisms $\tau_J$ and restrictions}

\subsection{Restriction to $\uqg_J$}
Given a subset $J$ of $I$, we denote by $\uqg_J$ the subalgebra of
$\uqg$ generated by $ x^\pm_{i,n},\wt{k}_i^{\pm 1}, h_{i,r}, i\in J, n
\in \Z, r\in \Z\backslash 0$. Let $$\on{res}_J: \on{Rep} \uqg \arr \on{Rep} \uqg_J$$ be the
restriction map and $\beta_J$ be the homomorphism $\Z[Y_{i,a}^{\pm
1}]_{i\in I; a \in \C^\times} \arr \Z[Y_{i,a}^{\pm 1}]_{i\in J; a \in
\C^\times}$, sending $Y_{i,a}^{\pm 1}$ to itself for $i \in J$ and to
$1$ for $i \not{\hspace*{-1mm}\in} J$.

According to Theorem 3(3) of \cite{FR:char}, the diagram
$$\begin{CD}
\on{Rep} \uqg @>{\chi_q}>>  \Z[Y_{i,a}^{\pm 1}]_{i\in I; a \in
\C^\times}\\
@VV{\on{res}_J}V     @VV{\beta_J}V\\
\on{Rep} \uqg_J  @>{\chi_{q,J}}>>  \Z[Y_{i,a}^{\pm 1}]_{i\in J; a \in
\C^\times}
\end{CD}$$
is commutative.

We will now refine the homomorphisms $\beta_J$ and $\on{res}_J$.

\subsection{The homomorphism $\tau_J$}    \label{homtauJ}

Consider the elements $\wt{h}_{i,n}$ defined by formula
\eqref{connection} and $\wt{k}_i^{\pm 1}$ defined by formula
\eqref{wtk}.

\begin{lem}
\begin{align*}
\wt{k}_i x^\pm_{j,n} \wt{k}_i^{-1} &= q^{\pm r_i \delta_{ij}}
x_{j,n}^{\pm},\\ [\wt{h}_{i,n},x_{j,m}^{\pm}] &= \pm \delta_{ij}
\frac{[nr_i]_{q}}{n} c^{\mp {|n|/2}} x_{j,n+m}^{\pm},\\
[\wt{h}_{i,n},h_{j,m}] &= \delta_{i,j}
\delta_{n,-m}\frac{[nr_i]_{q}}{n} \frac{c^n-c^{-n}}{q-q^{-1}}.
\end{align*}
In particular, $\wt{k}_i^{\pm 1}, \wt{h}_{i,n}, i \in \ol{J}, n \in
\Z\backslash 0$, where $\ol{J} = I - J$, commute with the subalgebra
$\uqg_J$ of $\uqg$.
\end{lem}

\begin{proof} These formulas follow from the relations
given in \thmref{defining} and the formula $B(q) \wt{C}(q) = D(q)$.
\end{proof}

Denote by $\uqhJ$ the subalgebra of $\uqg$ generated by $\wt{k}_i^{\pm
1}, \wt{h}_{i,n}, i \in \ol{J}, n \in \Z \setminus 0$. Then $\uqg_J
\otimes \uqhJ$ is naturally a subalgebra of $\uqg$. We can therefore
refine the restriction from $\uqg$--modules to $\uqg_J$--modules by
considering the restriction from $\uqg$--modules to $\uqg_J \otimes
\uqhJ$--modules.

Thus, we look at the common (generalized) eigenvalues of the operators
$k_i^{\pm 1}, h_{i,n}, i \in J$, and $\wt{k}_i^{\pm 1}, \wt{h}_{i,n},
i \in\ol{J}$. We know that the eigenvalues of $h_{i,n}$ have the form
\eqref{h eig}.  The corresponding eigenvalue of $\wt{h}_{i,n}$ equals
\begin{equation}    \label{wth eig}
\frac{[n]_q}{n} \sum_{j\in I} \wt{C}_{ji}(q^n) [r_j]_{q^n} \left(
\sum_{r=1}^{k_j} (a_{jr})^n - \sum_{s=1}^{l_j} (b_{js})^n \right),
\quad \quad n>0.
\end{equation}

According to \lemref{inverse lemma}, $\wt{C}_{ji}(x) =
\wt{C}'_{ji}(x)/d(x)$, where $\wt{C}'_{ji}(x)$ and $d(x)$ are certain
polynomials with positive integral coefficients (we fix a choice of
such $d(x)$ once and for all). Therefore formula \eqref{wth eig}
can be rewritten as
\begin{equation}    \label{wth eig1}
\frac{[n]_q}{nd(q^n)} \left( \sum_{m=1}^{u_i} (c_{im})^n -
\sum_{p=1}^{t_i} (d_{ip})^n \right),
\end{equation}
where $c_{im}$ and $d_{ip}$ are certain complex numbers (they are
obtained by multiplying $a_{jr}$ and $b_{js}$ with all monomials
appearing in $\wt{C}'_{ji}(q) [r_j]_{q}$).

According to \propref{common spectra}, to each monomial
\eqref{monomial} in $\chi_q(V)$ corresponds a generalized eigenspace
of $h_{i,n}, i\in I, n \in \Z \setminus 0$, with the common
eigenvalues given by formula \eqref{h eig} (note that the eigenvalues
of $k_i, i \in I$, can be read off from the weight of the
monomial). Using formula \eqref{wth eig} we find the corresponding
eigenvalues of $\wt{h}_{i,n}, i \in \ol{J}$ in the form \eqref{wth
eig1}. Now we attach to these common eigenvalues the following
monomial in the letters $Y_{i,a}^{\pm 1}, i \in J$, and $Z_{j,c}^{\pm
1}, j \in \ol{J}$:
$$
\left( \prod_{i\in J} \prod_{r=1}^{k_i} Y_{i,a_{ir}} \prod_{s=1}^{l_i}
Y_{i,b_{is}}^{-1} \right) \cdot \left( \prod_{k\in \ol{J}}
\prod_{m=1}^{u_k} Z_{k,c_{km}} \prod_{p=1}^{t_k} Z_{k,d_{kp}}^{-1}
\right).
$$

The above procedure can be interpreted as follows. Introduce the
notation
\bean    \label{yy}
\yy = \Z[Y_{i,a}^{\pm 1}]_{i\in I,a \in \C^\times},
\eean
\bean    \label{yyJ}
\yy^{(J)} = \Z[Y_{i,a}^{\pm 1}]_{i\in J,a \in
\C^\times} \otimes \Z[Z_{k,c}^{\pm 1}]_{k\in \ol{J},c\in
\C^\times}
\eean
Write
$$
(D(q) \wt{C}'(q))_{ij} = \sum_{k \in \Z} p_{ij}(k) q^k.
$$

\begin{definition}
The homomorphism $\tau_J: \yy \arr \yy^{(J)}$ is defined by the
formulas
\begin{align}    \label{tauJY1}
\tau_J(Y_{i,a}) &= Y_{i,a} \cdot \prod_{j \in \ol{J}} \prod_{k \in \Z}
Z^{p_{ij}(k)}_{j,aq^k}, \quad \quad i \in J,\\ \label{tauJY2}
\tau_J(Y_{i,a}) &= \prod_{j \in \ol{J}} \prod_{k \in \Z}
Z^{p_{ij}(k)}_{j,aq^k}, \quad \quad i \in \ol{J}.
\end{align}
\qed
\end{definition}

Observe that the homomorphism $\beta_J$ can be represented as the
composition of $\tau_J$ and the homomorphism $\yy^{(J)} \arr
\Z[Y_{i,a}^{\pm 1}]_{i\in J,a \in \C^\times}$ sending all $Z_{k,c}, k
\in \ol{J},$ to $1$. Therefore $\tau_J$ is indeed a refinement of
$\tau_J$, and so the restriction of $\tau_J$ to the image of $\on{Rep}
\uqg$ in $\yy$ is a refinement of the restriction homomorphism
$\on{res}_J$.

\subsection{Properties of $\tau_J$}

The main advantage of $\tau_J$ over $\beta_J$ is the following.

\begin{lem}    \label{tau injective1}
The homomorphism $\tau_J$ is injective.
\end{lem}

\begin{proof} The statement of the lemma follows from the fact that
the matrix $\wt{C}'(q)$ is non-degenerate.
\end{proof}

\begin{lem}    \label{restriction}
Let us write $\chi_q(V)$ as the sum $\sum_k P_k Q_k$, where $P_k \in
\Z[Y_{i,a}^{\pm 1}]_{i\in J,a \in \C^\times}$, $Q_k$ is a monomial in
$\Z[Z_{j,c}^{\pm 1}]_{j\in \ol{J},c\in \C^\times}$, and all monomials
$Q_k$ are distinct. Then the restriction of $V$ to $\uqg_J$ is
isomorphic to $\oplus_k V_k$, where $V_k$'s are $\uqg_J$--modules with
$\chi^J_q(V_k) = P_k$. In particular, there are no extensions between
different $V_k$'s in $V$.
\end{lem}

\begin{proof}
The monomials in $\chi_q(V) \in \yy$ encode the common eigenvalues of
$h_{i,n}, i \in I$ on $V$. It follows from \secref{homtauJ} that the
monomials in $\tau_J(\chi_q(V))$ encode the common eigenvalues of
$h_{i,n}, i \in J$, and $\wt{h}_{j,n}, j \in \ol{J}$, on $V$.

Therefore we obtain that the restriction of $V$ to $\uqg_J \otimes
\uqhJ$ has a filtration with the associated graded factors $V_k
\otimes W_k$, where $V_k$ is a $\uqg_J$--module with $\chi^J_q(V_k) =
P_k$, and $W_k$ is a one--dimensional $\uqhJ$--module, which corresponds
to $Q_k$. By our assumption, the modules $W_k$ over $\uqhJ$ are
pairwise distinct. Because $\uqhJ$ commutes with $\uqg_J$, there are
no extensions between $V_k \otimes W_k$ and $V_l \otimes W_l$ for
$k\neq l$, as $\uqg_J \otimes \uqhJ$--modules. Hence the restriction
of $V$ to $\uqg_J$ is isomorphic to $\oplus_k V_k$.
\end{proof}

Write
$$d(q)[r_i]_q = \sum_{k \in \Z} s_i(k) q^k.$$ Set
$$
B_{i,a} = \prod_{k \in \Z} Z_{i,aq^k}^{s_i(k)}.
$$

\begin{lem}    \label{tau injective2}
We have:
\begin{align}     \label{tauJA1}
\tau_J(A_{i,a}) &= \beta_J(A_{i,a}), \quad \quad i \in J,\\
\label{tauJA2}
\tau_J(A_{i,a}) &= \beta_J(A_{i,a}) B_{i,a}, \quad \quad i \in
\ol{J}
\end{align}
\end{lem}

\begin{proof} This follows from the formula $D(q) \wt{C}'(q) C(q) =
D(q) d(q)$.
\end{proof}

In the case when $J$ consists of a single element $j \in I$, we will
write $\yy^{(J)}, \tau_J$ and $\beta_J$ simply as $\yy^{(j)}, \tau_j$
and $\beta_j$. Consider the diagram (we use the notation \eqref{yy},
\eqref{yyJ}):

\bean\label{A-diagram} 
\renewcommand{\arraystretch}{1.5}
\begin{array}{ccc}
\yy & \stackrel{\tau_j}
{\longrightarrow} & \yy^{(j)} \\
{\downarrow }&&{\downarrow\;\; \ol{A}_{j,x}^{-1}} \\
\yy & \stackrel{\tau_j}
{\longrightarrow} & \yy^{(j)}
\end{array}
\eean where the map corresponding to the right vertical row is the
multiplication by \newline $\beta_j(A_{j,x})^{-1}\otimes 1$.

The following result will allow us to reduce various statements to the
case of $\uqsl$.

\begin{lem}\label{A}
There exists a unique map $\yy \arr \yy$, which makes the diagram
\Ref{A-diagram} commutative. This map is the multiplication by
$A^{-1}_{j,x}$.
\end{lem}

\begin{proof}
The fact that multiplication by $A^{-1}_{j,x}$ makes the diagram
commutative follows from formula \eqref{tauJA1}. The uniqueness
follows from the fact that $\tau_j$ and the multiplication by
$\beta_j(A_{j,x})^{-1}\otimes 1$ are injective maps.
\end{proof}

\section{The structure of $q$--characters}    \label{str}

In this section we prove Conjecture 1 from \cite{FR:char}.

Let $V$ be an irreducible finite-dimensional $\uqg$ module $V$
generated by highest weight vector $v$. Then by Proposition 3 in
\cite{FR:char}, \bean\label{form of char} \chi_q(V) = m_+( 1 +
\sum_p M_p), \eean where each $M_p$ is a monomial in $A_{i,c}^{\pm
1}$, $c\in\C^\times$ and $m_+$ is the highest weight monomial.

In what follows, by a monomial in $\Z[x_\al^{\pm 1}]_{\al \in A}$ we
will always understand a monomial in reduced form, i.e., one that does
not contain factors of the form $x_\al x_\al^{-1}$. Thus, in
particular, if we say that a monomial $M$ contains $x_\al$, it means
that there is a factor $x_\al$ in $M$ which can not be cancelled.

\begin{thm}\label{1} The $q$--character of an irreducible
finite-dimensional $\uqg$ module $V$ has the form \Ref{form of char}
where each $M_p$ is a monomial in $A_{i,c}^{-1}$, $i\in I$,
$c\in\C^\times$ (i.e., it does not contain any factors $A_{i,c}$).
\end{thm}

\begin{proof}
The proof follows from a combination of Lemmas \ref{tau injective1},
\ref{A} and \ref{inverse lemma}.

First, we observe that it suffices to prove the statement of Theorem
\ref{1} for fundamental representations $V_{\om_i}(a)$. Indeed, then
Theorem \ref{1} will be true for any tensor product of the fundamental
representations. By \corref{generated}, any irreducible representation
$V$ can be represented as a quotient of a submodule of a tensor
product $W$ of fundamental representations, which is generated by the
highest weight vector. Therefore each monomial in a $q$--character of
$V$ is also a monomial in the $q$--character of $W$. In addition, the
highest weight monomials of the $q$--characters of $V$ and $W$
coincide. This implies that Theorem \ref{1} holds for $V$.

Second, Theorem \ref{1} is true for $\g=\uqsl$. Indeed, by the
argument above, it suffices to check the statement for the
fundamental representation $V_1(a)$. But its $q$--character is known
explicitly (see \cite{FR:char}, formula (4.3)): \bean\label{sl2 fund}
\chi_q(V_1(a)) = Y_a+Y^{-1}_{aq^2}=Y_a(1+A_{aq}^{-1}).  \eean and it
satisfies the required property.

For general quantum affine algebra $\uqg$, we will prove Theorem
\ref{1} (for the case of the fundamental representations) by
contradiction.

Suppose that the theorem fails for some fundamental representation
$V_{\om_{i_0}}(a_0)=V$ and denote by $\chi$ its $q$--character
$\chi_q(V)$. Denote by $m_+$ the highest weight monomial $Y_{i_0,a}$
of $\chi$.

Recall from \secref{fd} that we have a partial order on the weight
lattice. It induces a partial order on the monomials occurring in
$\chi$. Let $m$ be the highest weight monomial in $\chi$, such that $m$
can not be written as a product of $m_+$ with a monomial in
$A_{i,c}^{-1}$, $i\in I$, $c\in\C^\times$. This means that
\begin{equation}    \label{assum}
\mbox{any monomial $m'$ in $\chi$, such that $m'>m$, is a product
of $m_+$ and $A_{i,c}^{-1}$'s.}
\end{equation}

In Lemmas \ref{Claim 1.} and \ref{Claim 2.} we will
establish certain properties of $m$ and in \lemref{Claim 3.} we
will prove that these properties can not be satisfied simultaneously.

Recall that a monomial in $\Z[Y_{i,a}^{\pm 1}]_{i\in I,a\in\C^\times}$
is called dominant if does not contain factors $Y_{i,a}^{-1}$ (i.e.,
if it is a product of $Y_{i,a}$'s in positive powers only).

\begin{lem}\label{Claim 1.}
The monomial $m$ is dominant.
\end{lem}

\begin{proof} Suppose $m$ is not dominant. Then it contains a factor
of the form $Y^{-1}_{i,a}$, for some $i\in I$. Consider
$\tau_i(\chi)$. By Lemma \ref{restriction}, we have \be \tau_i(\chi) =
\sum_p \chi_{q_i}(V_p) \cdot N_p, \ee where $V_p$'s are representation
of $\uqsli = \uqg_{\{ i \}}$ and $N_p$'s are monomials in
$Z^{\pm1}_{j,a}, j \neq i$. We have already shown that Theorem \ref{1}
holds for $\uqsli$, so \bean
\label{restriction improved} \tau_i(\chi)
=\sum_p \left( m_p(1+\sum_r \ol{M}_{r,p}) \right) \cdot N_p, \eean
where each $m_p$ is a product of $Y_{i,b}$'s (in positive powers
only), and each $\ol{M}_{r,p}$ is a product of several factors
$\ol{A}_{i,c}^{-1}=Y_{i,cq^{-1}}^{-1}Y_{i,cq}^{-1}$ (note that
$\ol{M}_{r,p} = \tau_i(M_{r,p})$.

Since $m$ contains $Y^{-1}_{i,a}$ by our assumption, the monomial
$\tau_i(m)$ is not among the monomials $\{m_p \cdot N_p\}$. Hence \be
\tau_i(m) = m_{p_0}\ol{M}_{r_0,p_0} \cdot N_{p_0}, \ee for some $p_0,
r_0$ and $\ol{M}_{r_0,p_0} \ne 1$. There exists a monomial $m'$ in
$\chi$, such that $\tau_i(m') = m_{p_0} \cdot N_{p_0}$. Therefore
using Lemma \ref{A} we obtain that \be m = m'M_{r_0,p_0}, \ee where
$M_{r_0,p_0}$ is obtained from $\ol M_{r_0,p_0}$ by replacing all
$\ol{A}^{-1}_{i,c}$ by ${A}_{i,c}^{-1}$. In particular, $m'>m$ and by
our assumption \eqref{assum} it can be written as $m'=m_+M'$, where
$M'$ is a product of $A^{-1}_{k,c}$. But then
$m=m'M_{r_0,p_0}=m_+M'M_{r_0,p_0}$, and so $m$ can be written as a
product of $m_+$ and a product of factors $A_{k,c}^{-1}$. This is a
contradiction. Therefore $m$ has to be dominant.
\end{proof}

\begin{lem}\label{Claim 2.} The monomial $m$ can be written in the form
\bean\label{one A} m=m_+ M \prod_pA_{j_0,a_p}, \eean where $M$ is a
product of factors $A_{i,c}^{-1}$, $i\in I$, $c\in\C^\times$. In other words,
if $m$ contains factors $A_{j,a}$, then all such $A_{j,a}$ have the
same index $j=j_0$.
\end{lem}

\begin{proof} Suppose that $m=m_+M$, where $M$ contains a factor
$A_{i,c}$. Let $V_m$ be the generalized eigenspace of the operators
$k_j^{\pm 1}, h_{j,n}, j \in I$, corresponding to the monomial $m$. We
claim that for all $v \in V_m$ we have: \bean\label{killed} x^+_{j,n}
\cdot v = 0, \qquad j\in I, j\neq i,\qquad n\in\Z.  \eean Indeed, let
$\tau_j(m) = \beta_j(m) \cdot N$ (recall that $\beta_j(m)$ is obtained
from $m$ by erasing all $Y_{s,c}$ with $s\ne j$ and $N$ is a monomial
in $Z^{\pm 1}_{s,c}$, $s\in I$, $s\ne j$).  By \lemref{restriction},
$x^+_{j,n} \cdot v$ belongs to the direct sum of the generalized
eigenspaces $V_{m'_p}$, corresponding to the monomials $m'_p$ in
$\chi$ such that $\tau_j(m'_p) = \beta_j(m'_p) \cdot N$ (with the same
$N$ as in $\tau_j(m) = \beta_j(m) \cdot N$). By formula
\eqref{tauJA2},
$$
\tau_j\left( m_+ \prod A_{i_k,c_k}^{\pm 1} \right) = \tau_j(m_+) \prod
\beta_j(A_{i_k,c_k})^{\pm 1} \prod_{i_k\neq j} B_{i_k,c_k}^{\pm 1}.
$$

In particular, $N$ contains a factor $B_{i,c}$, and therefore all
monomials $m'_p$ with the above property must contain a factor
$A_{i,c}$. By our assumption \eqref{assum}, the weight of each
$m'_p$ can not be higher then the weight of $m$. But the weight of
$x^+_{j,n} \cdot v$ should be greater than the weight of
$m$. Therefore we obtain formula \Ref{killed}.

Now, if $M$ contained factors $A_{i,c}$ and $A_{j,d}$ with $i\ne j$,
then any non-zero eigenvector (not generalized) in the generalized
eigenspace $V_m$ corresponding to $m$ would be a highest weight vector
(see formula \eqref{hwv}). Such vectors do not exist in $V$, because
$V$ is irreducible. The statement of the lemma now follows.
\end{proof}

\begin{lem}\label{Claim 3.} Let $m$ be any monomial in the $q$--character
of a fundamental representation that can be written in the form \Ref{one
A}. Then $m$ is not dominant.
\end{lem}

\begin{proof} We say a monomial $M\in\yy$ (see \eqref{yy}) has {\em
lattice support with base $a_0\in\C^\times$} if
$M\in\Z[Y_{i,a_0q^k}^{\pm1}]_{i\in I,k\in\Z}$.

Any monomial $m \in \yy$ can be uniquely written as a product
$m=m^{(1)}\dots m^{(s)}$, where each monomial $m^{(i)}$ has lattice
support with a base $a_i$, and $a_i/a_j\not\in q^\Z$ for $i\neq
j$. Note that a non-constant monomial in $A_{i,bq^k}^{\pm 1},i\in
I,k\in\Z$, can not be equal to a monomial in $A_{i,cq^k}^{\pm 1},i\in
I,k\in\Z$ if $b/c \not\in q^{\Z}$. Therefore if $m$ can be written in
the form \Ref{one A}, then each $m^{(i)}$ can be written in the form
\Ref{one A}, where $m_+ = Y_{i_0,a}$ if $a_i=a$, and $m_+ = 1$ if
$a/a_i\not\in q^\Z$ (note that the product over $p$ in \eqref{one A}
may be empty for some $m^{(i)}$). We will prove that none of
$m^{(i)}$'s is dominant unless $m^{(i)}=m_+$ or $m^{(i)}=1$.

Consider first the case of $m^{(1)}$, which has lattice support with
base $a$. Then \be m^{(1)} = \prod_{i\in I} \prod_{n\in
\Z}Y^{p_{i}(n)}_{i,aq^n}.  \ee Define Laurent polynomials $P_i(x)$,
$i\in I$ by \be P_i(x)=\sum_{n\in\Z} p_{i}(n) x^n.  \ee If $m^{(1)}$
can be written in the form \Ref{one A}, then \bean\label{polynomials}
P_i(x)=-\sum_{j\in I}C_{ij}(x)R_j(x)+\delta_{i,i_0},\qquad \forall
i\in I, \eean where $R_j(x)$'s are some polynomials with integral
coefficients. All of these coefficients are non-negative if $j\neq
j_0$. Now suppose that $m^{(1)}$ is a dominant monomial.  Then each
$P_i(x)$ is a polynomial with non-negative coefficients. We claim that
this is possible only if all $R_i(x) = 0$.

Indeed, according to Lemma \ref{inverse lemma}, the coefficients of
the inverse matrix to $C(x)$, $\wt C(x)$, can be written in the form
\Ref{inverse formula}, where $\wt C_{jk}'(x)$, $d(x)$ are polynomials
with non-negative coefficients. Multiplying \Ref{polynomials} by $\wt
C'(x)$, we obtain \bean \label{pol1} \sum_{j\in I} P_j(x) \wt
C_{jk}'(x)+d(x)R_k(x)=\wt C_{i_0,k}'(x), \qquad \forall k\in I. \eean

Given a Laurent polynomial
$$
p(x) = \sum_{-r\leq i\leq s} p_i x^i, \quad \quad p_{-r} \neq 0,
p_s\neq 0,
$$
we will say that the length of $p(x)$ equals $r+s$. Clearly, the
length of the sum and of the product of two polynomials with
non-negative coefficients is greater than or equal to the length of
each of them. Therefore if $k\neq j_0$, and if $R_k(x)\neq 0$, then
the length of the LHS is greater than or equal to the length of
$d(x)$, which is greater than the length of $\wt C_{i_0,k}'$ by Lemma
\ref{inverse lemma}. This implies that $R_k(x)=0$ for $k\neq j_0$.

Hence $m^{(1)}$ can be written in the form \be m^{(1)} = Y_{i,a}
\prod_{n\in \Z} A^{c_n}_{j_0,aq^n}.  \ee But such a monomial can not be
dominant because its weight is $\omega_i - n \al_{j_0}$, where $n>0$,
and such a weight is not dominant. This proves the required statement
for the factor $m^{(1)}$ of $m$ (which has lattice support with base
$a$).

Now consider a factor $m^{(i)}$ with lattice support with base $b$,
such that $b/a\not\in q^\Z$. In this case we obtain the following
equation: the LHS of formula \eqref{pol1} $= 0$. The previous
discussion immediately implies that there are no solutions of this
equation with non-zero polynomials $R_k(x)$ satisfying the above
conditions. This completes the proof of the lemma.
\end{proof}

Theorem \ref{1} now follows from Lemmas \ref{Claim 1.}, \ref{Claim 2.}
and \ref{Claim 3.}.
\end{proof}

\begin{cor}\label{no dominant monomials}
The only dominant monomial in $\chi_q(V_{\omega_i}(a))$ is the highest
weight monomial $Y_{i,a}$.
\end{cor}
\begin{proof} 
This follows from the proof of \lemref{Claim 3.}.
\end{proof}

\section{A characterization of $q$--characters in terms of the
screening operators}

In this section we prove Conjecture 2 from \cite{FR:char}.

\subsection{Definition of the screening operators}    \label{defscr}

First we recall the definition of the screening operators on $\yy =
\Z[Y_{i,a}^{\pm 1}]_{i\in I; a \in \C^\times}$ from \cite{FR:char} and
state the main result.

Consider the free $\yy$--module with generators $S_{i,x}, x \in
\C^\times$,
$$
\wt{\yy}_i = \us{x \in \C^\times}{\opl} \yy \cdot S_{i,x}.
$$
Let $\yy_i$ be the quotient of $\wt{\yy}_i$ by the relations
\begin{equation}    \label{relation}
S_{i,xq_i^2} = A_{i,xq_i} S_{i,x}.
\end{equation}
Clearly,
$$
\yy_i \simeq \us{x \in (\C^\times/q_i^{2\Z})}{\opl} \yy \cdot S_{i,x},
$$
and so $\yy_i$ is also a free $\yy$--module.

Define a linear operator $\wt{S}_i: \yy \arr \wt{\yy}_i$ by the
formula
$$
\wt{S}_i(Y_{j,a}) = \delta_{ij} Y_{i,a} S_{i,a}
$$
and the Leibniz rule: $\wt{S}_i(ab) = b \wt{S}_i(a) + a
\wt{S}_i(b)$. In particular,
$$
\wt{S}_i(Y_{j,a}^{-1}) = - \delta_{ij} Y_{i,a}^{-1} S_{i,a}.
$$

Finally, let $$S_i: \yy \arr \yy_i$$ be the composition of $\wt{S}_i$
and the projection $\wt{\yy}_i \arr \yy_i$. We call $S_i$ the $i$th
{\em screening operator}.

The following statement was conjectured in \cite{FR:char} (Conjecture
2).

\begin{thm}    \label{main}
The image of the homomorphism $\chi_q$ equals the intersection of
the kernels of the operators $S_i, i \in I$.
\end{thm}

In \cite{FR:char} this theorem was proved in the case of $\uqsl$.
In the rest of this section we prove it for an arbitrary $\uqg$.

\subsection{Description of $\on{Ker} S_i$}

First, we describe the kernel of $S_i$ on $\yy$. The following result
was announced in \cite{FR:char}, Proposition 6.

\begin{prop}    \label{keri}
The kernel of $S_i: \yy \arr \yy_i$ equals 
\bean\label{Ri}
{\mc K}_i =
\Z[Y_{j,a}^{\pm 1}]_{j\neq i; a \in \C^\times} \otimes \Z[Y_{i,b}
+ Y_{i,b} A_{i,bq_i}^{-1}]_{b \in \C^\times}.
\eean
\end{prop}

\begin{proof}
A simple computation shows that ${\mc K}_i \subset \on{Ker}_{\yy}
S_i$. Let us show that $\on{Ker}_{\yy} S_i \subset {\mc K}_i$.

For $x \in \C^\times$, denote by $\yy(x)$ the subring
$\Z[Y_{j,xq^{n}}^{\pm 1}]_{j \in I, n \in \Z}$ of $\yy$. We have:
$$
\yy \simeq \us{x \in (\C^\times/q^{\Z})}{\otimes} \yy(x).
$$

\begin{lem}    \label{split}
$$
\on{Ker}_{\yy} S_i = \us{x \in (\C^\times/q^{\Z})}{\otimes}
\on{Ker}_{\yy(x)} S_i.
$$
\end{lem}

\begin{proof} Let $P \in \yy$, and
suppose it contains $Y_{j,a}^{\pm 1}$ for some $a \in \C^\times$ and
$j \in I$. Then we can write $P$ as the sum $\sum_k R_k Q_k$, where
$Q_k$'s are distinct monomials, which are products of the factors
$Y_{s,aq^n}^{\pm 1}, s \in I, n \in \Z$ (in particular, one of the
$Q_k$'s could be equal to $1$), and $R_k$'s are polynomials which do
not contain $Y_{s,aq^n}^{\pm 1}, s \in I, n \in \Z$. Then
$$
S_i(P) = \sum_k (Q_k \cdot S_i(R_k) + R_k \cdot S_i(Q_k)).
$$
By definition of $S_i$, $S_i(Q_k)$ belongs to $\yy \cdot S_{i,a}$,
while $S_i(R_k)$ belongs to the direct sum of $\yy \cdot S_{i,b}$,
where $b \not\in aq^{\Z}$.

Therefore if $P \in \on{Ker}_{\yy} S_i$, then $\sum_k Q_k \cdot
S_i(R_k) = 0$. Since $Q_k$'s are distinct, we obtain that $R_k \in
\on{Ker}_{\yy} S_i$. But then $S_i(P) = \sum_k R_k \cdot
S_k(Q_k)$. Therefore $P$ can be written as $\sum_l R_l \wt{Q}_l$,
where each $\wt{Q}_l$ is a linear combination of the $Q_k$'s, such
that $\wt{Q}_l \in \on{Ker}_{\yy} S_i$. This proves that
$$P \in \on{Ker}_{\yy(\neq a)} S_i \otimes \on{Ker}_{\yy^{(a)}} S_i,$$
where $\yy(\neq \hspace*{-2mm} a) = \Z[Y_{j,b}^{\pm 1}]_{j \in I, b
\not\in aq^{\Z}}$. By repeating this procedure we obtain the lemma
(because each polynomial contains a finite number of variables
$Y_{j,a}^{\pm 1}$, we need to apply this procedure finitely many
times).
\end{proof}

According to \lemref{split}, it suffices to show that
$\on{Ker}_{\yy(x)} S_i \subset {\mc K}_i(x)$, where
$$
{\mc K}_i(x) = \Z[Y_{j,xq^n}^{\pm 1}]_{j\neq i; n \in \Z} \otimes
\Z[Y_{i,xq^n} + Y_{i,xq^n} A_{i,xq^{n}q_i}^{-1}]_{n \in \Z}.
$$
Denote $Y_{j,xq^n}$ by $y_{j,n}$, $A_{j,xq^n}$ by $a_{j,n}$, and
$A_{j,xq^n} Y_{j,xq^{n}q_j}^{-1} Y_{j,xq^{n}q_j^{-1}}^{-1}$ by
$\ol{a}_{j,n}$. Note that $\ol{a}_{j,n}$ does not contain factors
$y_{j,m}^{\pm 1}, m \in \Z$.

Let $T$ be the shift operator on $\yy(x)$ sending $y_{j,n}$ to
$y_{j,n+1}$ for all $j \in I$. It follows from the definition of $S_i$
that $P \in \on{Ker}_{\yy(x)} S_i$ if and only if $T(P) \in
\on{Ker}_{\yy(x)} S_i$. Therefore (applying $T^m$ with large enough
$m$ to $P$) we can assume without loss of generality that $P \in
\Z[y_{i,n},y_{i,n+2r_i}^{-1}]_{n\geq 0} \otimes \Z[y_{j,n}^{\pm
1}]_{j\neq i, n\geq 0}$.

We find from the definition of $S_i$:
\begin{align} \notag
S_i(y_{j,n}) &= 0, \quad \quad j \neq i,\\ \label{appl}
S_i(y_{i,2r_i n+\ep}) &= y_{i,\ep} \prod_{k=1}^n y_{i,2r_i k+\ep}^2
\ol{a}_{i,r_i(2k-1)+\ep} \cdot S_{i,xq^\ep},
\end{align}
where $\ep \in \{ 0,1,\ldots,2r_i-1\}$. Therefore each $P \in
\on{Ker}_{\yy(x)} S_i$ can be written as a sum $P = \sum P_\ep$,
where each $P_\ep \in \on{Ker}_{\yy(x)} S_i$ and
$$
P_\ep \in \Z[y_{i,2r_i n+\ep},y_{i,2r_i(n+1)+\ep}^{-1}]_{n\geq 0} \otimes
\Z[y_{j,n}^{\pm 1}]_{j\neq i, n\geq 0}.
$$

It suffices to consider the case $\ep=0$. Thus, we show that if
$$
P \in \yy_i^{\geq 0}(x) = \Z[y_{i,2r_i n},y_{i,2r_i(n+1)}^{-1}]_{n\geq
0} \otimes \Z[y_{j,n}^{\pm 1}]_{j\neq i, n\geq 0},
$$
then
$$
P \in {\mc K}_i^{\geq 0}(x) = \Z[t_{n}]_{n\geq 0} \otimes
\Z[y_{j,n}^{\pm 1}]_{j\neq i, n\geq 0},
$$
where
$$
t_{n} = y_{i,2r_i n} + y_{i,2r_i n} a_{i,r_i(2n+1)}^{-1} = y_{i,2r_i
n} + y_{i,2r_i(n+1)}^{-1} \ol{a}_{i,r_i(2n+1)}^{-1}.
$$

Consider a homomorphism ${\mc K}^{\geq 0}_{i}(x) \otimes \Z[y_{i,2r_i
n}]_{n\geq 0} \to \yy_i^{\geq 0}(x)$ sending $y_{j,n}^{\pm 1}, j \neq
i$ to $y_{j,n}^{\pm 1}$, $y_{i,2r_i n}$ to $y_{i,2r_i n}$, and $t_n$
to $y_{i,2r_i n} + y_{i,2r_i(n+1)}^{-1}
\ol{a}_{i,r_i(2n+1)}^{-1}$. This homomorphism is surjective, and its
kernel is generated by the elements
\begin{equation}    \label{Kernel}
(t_n - y_{i,2r_i n}) \ol{a}_{i,r_i(2n+1)} y_{i,2r_i(n+1)} - 1.
\end{equation}
Therefore we identify $\yy_i^{\geq 0}(x)$ with the quotient of ${\mc
K}^{\geq 0}_{i}(x) \otimes \Z[y_{i,2r_i n}]_{n\geq 0}$ by the ideal
generated by elements of the form \eqref{Kernel}.

Consider the set of monomials
$$
t_{n_1} \ldots t_{n_k} y_{i,2r_i m_1} \ldots y_{i,2r_i m_l}
\prod_{j\neq i,p_j\geq 0} y_{j,p_j}^{\pm 1},
$$
where all $n_1 \geq n_2 \geq \ldots n_k \geq 0, m_1 \geq m_2 \geq
\ldots m_l \geq 0$, and also $m_j \neq n_i+1$ for all $i$ and $j$. We
call these monomials {\em reduced}. It is easy to see that the set of
reduced monomials is a basis of $\yy_i^{\geq 0}(x)$.

Now let $P$ be an element of the kernel of $S_i$ on $\yy_i^{\geq 0}(x)$. Let us write it as a linear combination of the reduced
monomials. We represent $P$ as $y_{i,2r_i N}^a Q + R$. Here $N$ is the
largest integer, such that $y_{i,2r_i N}$ is present in at least one
of the basis monomials appearing in its decomposition; $a>0$ is the
largest power of $y_{i,2r_i N}$ in $P$; $Q \neq 0$ does not contain
$y_{i,2r_i N}$, and $R$ is not divisible by $y_{i,2r_i N}^a$. Recall
that here both $y_{i,2r_i N}^a Q$ and $R$ are linear combinations of
reduced monomials.

Recall that $S_i(t_n) = 0$, $S_i(y_{j,n}^{\pm 1}) = 0, j \neq i$,
and $S_i(y_{i,2r_i n})$ is given by formula \eqref{appl}.

Suppose that $N>0$. According to formula \eqref{appl},
\begin{equation}    \label{leading}
S_i(P) = a y_{i,2r_i N}^{a+1} \prod_{k=1}^{N-1} y_{i,2r_i k}
\prod_{l=1}^{N} \ol{a}_{i,r_i(2l-1)} y_{i,0} \; Q \cdot S_{i,x} +
\ldots
\end{equation}
where the dots represent the sum of terms that are not divisible by
$y_{i,2r_i N}^{a+1}$. Note that the first term in \eqref{leading} is
non-zero because the ring $\yy_i^{\geq 0}(x)$ has no divisors of
zero.

The monomials appearing in \eqref{leading} are not necessarily
reduced. However, by construction, $Q$ does not contain $t_{N-1}$, for
otherwise $y_{i,2r_i N}^a Q$ would not be a linear combination of
reduced monomials. Therefore when we rewrite \eqref{leading} as a
linear combination of reduced monomials, each reduced monomial
occurring in this linear combination is still divisible by $y_{i,2r_i
N}^{a+1}$. On the other hand, no reduced monomials occurring in the
other terms of $S_i(P)$ (represented by dots) are divisible by
$y_{i,2r_i N}^{a+1}$. Hence for $P$ to be in the kernel, the first
term of \eqref{leading} has to vanish, which is impossible. Therefore
$P$ does not contain $y_{i,2r_i m}$'s with $m>0$.

But then $P = \sum_k y_{i,0}^{p_k} R_k$, where $R_k \in {\mc K}_{i\geq
0}(x)$, and $S_i(P) = \sum_k p_k y_{i,0}^{p_k-1} R_k \cdot
S_{i,x}$. Such $P$ is in the kernel of $S_i$ if and only if all
$p_k=0$ and so $P \in {\mc K}_{i\geq 0}(x)$. This completes the
proof of \propref{keri}.
\end{proof}

Set \bean\label{Q} {\mc K}=\bigcap_{i\in I} {\mc K}_i=\bigcap_{i\in I}
\left(\Z[Y_{j,a}^{\pm 1}]_{j\neq i; a \in \C^\times} \otimes
\Z[Y_{i,b} + Y_{i,b} A_{i,bq_i}^{-1}]_{b \in \C^\times}\right).  \eean
Now we will prove that the image of the $q$--character homomorphism
$\chi_q$ equals ${\mc K}$.

\subsection{The image of $\chi_q$ is a subspace of ${\mc K}$}

First we show that the image of $\on{Rep} \uqg$ in $\yy$ under the
$q$--character homomorphism belongs to the kernel of $S_i$.

Recall the ring $\yy^{(i)} = \Z[Y_{i,a}^{\pm 1}]_{a \in \C^\times}
\otimes \Z[Z_{j,c}^{\pm 1}]_{j\neq i,c\in \C^\times}$ and the
homomorphism $\tau_i: \yy \arr \yy^{(i)}$ from \secref{homtauJ}.

Let $\ol{\yy}_i$ be the quotient of $\ds \us{x \in \C^\times}{\opl}
\Z[Y_{i,a}^{\pm 1}]_{a \in \C^\times} \cdot S_{i,x}$ by the submodule
generated by the elements of the form $S_{i,xq_i^2} - \ol{A}_{i,xq_i}
S_{i,x}$, where $\ol{A}_{i,xq_i} = Y_{i,x} Y_{i,xq_i^2}$.  Define a
derivation $\ol{S}_i: \Z[Y_{i,a}^{\pm 1}]_{a \in \C^\times} \arr
\ol{\yy}_i$ by the formula $\ol{S}_i(Y_{i,a}) = Y_{i,a}
S_{i,a}$. Thus, $\ol{\yy}_i$ coincides with the module $\yy_i$ in the
case of $\uqsli$ and $\ol{S}_i$ is the corresponding screening
operator.

Set $$\yy^{(i)}_i = \Z[Z_{j,c}^{\pm 1}]_{j\neq i,c\in \C^\times}
\otimes \ol{\yy}_i.$$ The map $\ol{S}_i$ can be extended uniquely to a
map $\yy^{(i)} \arr \yy^{(i)}_i$ by $\ol{S}_i(Z_{j,c}) = 0$
for all $j\neq i,c\in \C^\times$ and the Leibniz rule. We will also
denote it by $\ol{S}_i$. The embedding $\tau_i$ gives rise to an
embedding $\yy_i \arr \yy^{(i)}_i$ which we also denote by $\tau_i$.

\begin{lem}    \label{comm1}
The following diagram is commutative
$$
\begin{CD}
\yy @>{\tau_i}>> \yy^{(i)} \\
@VV{S_i}V     @VV{\ol{S}_i}V\\
\yy_i @>{\tau_i}>> \yy^{(i)}_i
\end{CD}
$$
\end{lem}

\begin{proof}
Since $\tau_i$ is a ring homomorphism and both $S_i$, $\ol{S}_i$ are
derivations, it suffices to check commutativity on the generators.
Let us choose a representative $x$ in each $q_i^{2\Z}$--coset of
$\C^\times$. Then we can write:
$$
\yy_i = \us{x \in \C^\times/q_i^{2\Z}}{\opl} \yy \cdot
S_{i,x}, \quad \quad \yy_i^{(i)} = \us{x \in
\C^\times/q_i^{2\Z}}{\opl} \yy^{(i)} \cdot S_{i,x}.
$$
By definition,
\begin{align*}
S_i(Y_{j,xq_i^{2n}}) &= \delta_{ij} Y_{i,x} \prod_m
A_{i,xq_i^{2m+1}}^{\pm 1} \; S_{i,x}, \\
\ol{S}_i(Y_{i,xq_i^{2n}}) &= Y_{i,x} \prod_m
\ol{A}_{i,xq_i^{2m+1}}^{\pm 1} \; S_{i,x}, \\
\ol{S}_i(Z_{j,c}) &= 0 , \quad \quad \forall j \neq i.
\end{align*}
Recall from formula \eqref{tauJY1} that $\tau_i(Y_{i,x})$ equals
$Y_{i,x}$ times a monomial in $Z_{j,c}^{\pm 1}, j\neq i$, and from
formula \eqref{tauJA2} that $\tau_i(A_{i,b}^{\pm 1}) =
\ol{A}_{i,b}^{\pm 1}$. Using these formulas we obtain:
$$
(\tau_i \circ S_i)(Y_{i,xq_i^{2n}}) = (\ol{S}_i \circ
\tau_i)(Y_{i,xq^{2n}}) = \tau_i(Y_{i,x}) \prod
\ol{A}_{i,xq_i^{2m+1}}^{\pm 1} \; S_{i,x}.
$$
On the other hand, when $j \neq i$, $\tau_i(Y_{j,x})$ is a monomial in
$Z_{k,c}^{\pm 1}, k \neq i$, according to formula
\eqref{tauJY2}. Therefore
$$
(\tau_i \circ S_i)(Y_{j,x}) = (\ol{S}_i \circ \tau_i)(Y_{j,x}) = 0,
\quad \quad j\neq i.
$$
The proves the lemma.
\end{proof}

\begin{cor}    \label{image}
The image of the $q$--character homomorphism $\chi_q: \on{Rep} \uqg
\arr \yy$ is contained in the kernel of $S_i$ on $\yy$.
\end{cor}

\begin{proof} Let $V$ be a finite-dimensional representation of
$\uqg$. We need to show that $S_i(\chi_q(V)) = 0$. By
\lemref{restriction}, we can write $\chi_q(V)$ as the sum $\sum_k P_k
Q_k$, where each $P_k \in \Z[Y_{i,a}^{\pm 1}]_{a\in \C^\times}$ is in
the image of the homomorphism $\chi_{q}^{(i)}: \on{Rep} \uqsli \arr
\Z[Y_{i,a}^{\pm 1}]_{a\in \C^\times}$, and $Q_k$ is a monomial in
$Z_{j,c}^{\pm 1}, j\neq i$.

The image of $\chi_{q}^{(i)}$ lies in the kernel of the operator
$\ol{S}_i$ (in fact, they are equal, but we will not use this
now). This immediately follows from the fact that $\on{Rep} \uqsl
\simeq \Z[\chi_q(V_1(a))]$ and $\ol{S}_i (\chi_q(V_1(a))) = 0$, which
is obtained by a straightforward calculation. We also have:
$\ol{S}_i(Z_{j,c}) = 0, \forall j \neq i$. Therefore $(\ol{S}_i \circ
\tau_i)(\chi_q(V)) = 0$. By \lemref{comm1}, $(\tau_i \circ
S_i)(\chi_q(V)) = 0$. Since $\tau_i$ is injective by \lemref{tau
injective1}, we obtain: $S_i(\chi_q(V)) = 0$.
\end{proof}

\subsection{${\mc K}$ is a subspace of the image of $\chi_q$}

Let $P \in {\mc K}$. We want to show that $P \in \on{Im} \chi_q$.

A monomial $m$ contained in $P \in \yy$ is called {\em highest
monomial} (resp., {\em lowest monomial}), if its weight is not lower
(resp., not higher) than the weight of any other monomial contained in
$P$.

\begin{lem}    \label{dom}
Let $P \in {\mc K}$. Then any highest monomial in $P$ is dominant and
any lowest weight monomial in $P$ is antidominant.
\end{lem}
\begin{proof}
First we prove that the highest monomials are dominant.

By \propref{keri}, $$P \in {\mc K}_i = \Z[Y_{j,a}^{\pm 1}]_{j\neq i; a
\in \C^\times} \otimes \Z[Y_{i,b} + Y_{i,b} A_{i,bq_i}^{-1}]_{b \in
\C^\times}.$$ The statement of the lemma will follow if we show that a
highest weight monomial contained in any element of ${\mc K}_i$ does not
contain factors $Y_{i,a}^{-1}$.

Indeed, the weight of $Y_{i,a}$ is $\omega_i$, and the weight of
$Y_{i,b} A_{i,bq_i}^{-1}$ is $\omega_i-\al_i$. Denote $t_b =
\Z[Y_{i,b} + Y_{i,b} A_{i,bq_i}^{-1}]_{b \in \C^\times}$. Given a
polynomial $Q \in \Z[t_b]_{b \in \C^\times}$, let $m_1,\ldots,m_k$ be
its monomials (in $t_b$) of highest degree. Clearly, the monomials of
highest weight in $Q$ (considered as a polynomial in $Y_{j,a}^{\pm
1}$) are $m_1,\ldots,m_k$, in which we substitute each $t_b$ by
$Y_{i,b}$. These monomials do not contain factors $Y_{i,a}^{-1}$.

The statement about the lowest weight monomials is proved similarly,
once we observe that $${\mc K}_i = \Z[Y_{j,a}^{\pm 1}]_{j\neq i; a \in
\C^\times} \otimes \Z[Y_{i,b}^{-1} + Y_{i,bq_i^{-2}}
A_{i,bq_i^{-1}}]_{b \in \C^\times}.$$
\end{proof}

Let $m$ be a highest monomial in $P$, and suppose that it enters $P$
with the coefficient $\nu_m \in \Z \setminus 0$.  Then $m$ is dominant
by \lemref{keri}. According to \thmref{ChP2}(2) and formula
\eqref{mp}, there exists an irreducible representation $V_1$ of
$\uqg$, such that $m$ is the highest weight monomial in
$\chi_q(V_1)$. Since $\chi_q(V_1) \in {\mc K}$ by \corref{image}, we
obtain that $P_1 = P - \nu_m \cdot \chi_q(V_1) \in {\mc K}$.

For $P \in \yy$, denote by $\Lambda(P)$ the (finite) set of dominant
weights $\la$, such that $P$ contains a monomial of weight greater
than or equal to $\la$. By \propref{keri}, if $P \in {\mc K}$ and
$\Lambda(P)$ is empty, then $P$ is necessarily equal to $0$.

Note that for any irreducible representation $V$ of $\uqg$ of highest
weight $\mu$, $\Lambda (\chi_q(V))$ is the set of all dominant weights
which are less than or equal to $\mu$. Therefore $\Lambda(P_1)$ is
properly contained in $\Lambda(P)$. By applying the above subtraction
procedure finitely many times, we obtain an element $P_k = \ds P -
\sum_{i=1}^k \chi_q(V_i)$, for which $\Lambda(P_k)$ is empty. But then
$P_k = 0$.

This shows that ${\mc K} \subset \on{Im} \chi_q$. Together with
\lemref{image}, this gives us \thmref{main} and the following
corollary.

\begin{cor}    \label{subring}
The $q$--character homomorphism,
\be
\chi_q:\; \on{Rep} \uqg \to {\mc K},
\ee
where ${\mc K}$ is given by \Ref{Q}, is a ring isomorphism.
\end{cor}

\subsection{Application: Algorithm for constructing
$q$-characters}\label{algorithm} 

Consider the following problem: Give an algorithm which for any
dominant monomial $m_+$ constructs the $q$--character of the
irreducible $\uqg$--module whose highest weight monomial is $m_+$.  In
this section we propose such an algorithm.  We prove that our
algorithm produces the $q$--characters of the fundamental
representations (in this case $m_+=Y_{i,a}$). We conjecture that
the algorithm works for any irreducible module.

Roughly speaking, in our algorithm we start from $m_+$ and gradually
expand it in all possible $U_{q_i}\widehat{\mathfrak {sl}}_2$
directions.  (Here we use the explicit formulas for $q$--characters of
$U_{q}\widehat{\mathfrak {sl}}_2$ and Lemma \ref{A}.)  In the process
of expansion some monomials may come from different directions. We
identify them in the maximal possible way.

First we introduce some terminology.

Let $\chi\in\Z_{\geq 0}[Y_{i,a}^{\pm1}]_{i\in I,a\in\C^\times}$ be a
polynomial and $m$  a monomial in
$\chi$ occuring with coefficient $s\in\Z_{>0}$. By definition, a
{\em coloring} of 
$m$ is a set $\{s_i\}_{i\in I}$ of non-negative integers such
that $s_i\leq s$. A polynomial $\chi$ in which all monomials are
colored is called a {\em colored polynomial}.  

We think of $s_i$ as the number of monomials of type $m$ which have
come from direction $i$ (or by expanding with respect to the $i$-th
subalgebra $U_{q_i}\widehat{\mathfrak {sl}}_2$).

A monomial $m$ is called $i$--{\em dominant} if it does not contain
variables $Y_{i,a}^{-1}$, $a\in\C^\times$. 

A monomial $m$ occurring in a colored polynomial $\chi$ with
coefficient $s$ is called {\em admissible} if
$m$ is $j$--dominant for all $j$ such that $s_j<s$. A colored
polynomial is called 
admissible if all of its monomials are admissible.

Given an admissible 
monomial $m$ occurring with coefficient $s$ in a colored polynomial
$\chi$, we define a new 
colored polynomial $i_m(\chi)$, called the $i$--expansion of $\chi$
with respect to $m$, as follows. 

If $s_i=s$, then $i_m(\chi)=\chi$. Suppose that $s_i<s$ and let $\ol
m$ be obtained from $m$ by setting $Y_{j,a}^{\pm1}=1$, for all $j\neq
i$. Since $m$ is admissible, $\ol m$ is a dominant monomial.
Therefore there exists an irreducible $U_{q_i}\widehat{\mathfrak
{sl}}_2$ module $V$, such that the highest weight monomial of $V$ is
$\ol m$. We have explicit formulas for the $q$-characters of all
irreducible $U_{q}\widehat{\mathfrak {sl}}_2$--modules (see, e.g.,
\cite{FR:char}, Section 4.1).  We write $\chi_{q_i}(V)=\ol m(1+\sum_p
{\ol M}_p)$, where $\ol M_p$ is a product of $\ol A_{i,a}^{-1}$. Let
\bean\label{mu form} \mu=m(1+\sum_p M_p), \eean where $M_p$ is
obtained from $\ol M_p$ by replacing all $\ol A_{i,a}^{-1}$ by
$A_{i,a}^{-1}$.

The colored polynomial $i_m(\chi)$ is obtained from $\chi$ by adding
monomials occuring in $\mu$ by the following rule. Let monomial $n$
occur in $\mu$ with coefficient $t\in\Z_{>0}$. If $n$ does not occur
in $\chi$ then it is added with the coefficient $t(s-s_i)$ and we set
the $i$-th coloring of $n$ to be $t(s-s_i)$, and the other colorings
to be $0$. If $n$ occurs in $\chi$ with coefficient $r$ and coloring
$\{r_i\}_{i\in I}$, then the new coefficient of $n$ in $i_m(\chi)$ is
$\max\{r,r_i+t(s-s_i)\}$. In this case the $i$-th coloring is changed
to $r_i+t(s-s_i)$ and other colorings are not changed.

Obviously, the $i$--expansions of $\chi$ with respect to $m$ commute
for different $i$. To expand a
monomial $m$ in all directions means to compute
$\ell_m(\dots2_m(1_m(\chi))\dots)$, where $\ell=rk(\g)$.

Now we describe the algorithm. We start with the
colored polynomial $m_+$ with all colorings set equal zero. Let the
$\uqgg$--weight of $m_+$ be $\la$. 
The set of weights of the form $\la-\sum_ia_i\al_i$, $a_i\in\Z_{\geq
0}$ has a natural partial order. Choose any total order compatible
with this partial order, so 
we have $\la=\la_1>\la_2>\la_3>\ldots$ 

At the first step we expand $m_+$ in all directions. Then we expand in
all directions all monomials of weight $\la_1$ obtained at the first
step.  Then we expand in all directions all monomials of weight
$\la_2$ obtained at the previous steps, and so on. Since the monomials
obtained in the expansion of a monomial of $\uqgg$--weight $\mu$ have
weights less than $\mu$, the result does not depend on the choice of
the total order.

Note that for any monomial $m$ except for $m_+$ occurring with
coefficient $s$ at any step, we have $\max_i\{s_i\}=s$. This property
means that we identify the monomials coming from different directions
in the maximal possible way.

The algorithm stops if all monomials have been expanded. 
We say that the algorithm fails at a monomial $m$
if $m$ is the first non-admissible monomial to be expanded.

Let $m_+$ be a dominant monomial and $V$ the corresponding irreducible
module.

\begin{conj}\label{alg conj}
The algorithm never fails and stops after finitely many
steps. Moreover, the final result of the algorithm is the
$q$--character of $V$.
\end{conj}

\begin{thm}
Suppose that  $\chi_q(V)$ does not contain dominant monomials other
then $m_+$. Then 
Conjecture \ref{alg conj} is true. In particular, Conjecture \ref{alg
conj} is true in the case of fundamental representations.
\end{thm}

\begin{proof}
For $i\in I$, let $D_i$ be a decomposition of the set of monomials in
$\chi_q(V)$ with multiplicities into a disjoint union of subsets such
that each subset forms the $q$--character of an irreducible
$U_{q_i}\widehat{\mathfrak {sl}}_2$ module. We refer to this
decompostion $D_i$ as the $i$-th decomposition of $\chi_q(V)$. Denote
$D$ the collection of $D_i$, $i\in I$.

Consider the following colored oriented graph $\Omega_V(D)$.  The
vertices are monomials in $\chi_q(V)$ with multiplicities.  We draw an
arrow of color $i$ from a monomial $m_1$ to a monomial $m_2$ if and
only if $m_1$ and $m_2$ are in the same subset of the $i$-th
decomposition and $m_2=A^{-1}_{i,a}m_1$ for some $a\in\C^\times$.

We call an oriented graph a tree (with one root) if there exists a
vertex $v$ (called root), such that there is an oriented path from $v$
to any other vertex.  The graph $\Omega_W(D)$, where $W$ is an
irreducible $U_{q}\widehat{\mathfrak {sl}}_2$--module is always a tree and
its root corresponds to the highest weight monomial.

Consider the full subgraph of $\Omega_V(D)$ whose vertices correspond
to monomials from a given subset of the $i$-th decomposition of
$\chi_q(V)$. All arrows of this subgraph are of color $i$. By Lemma
\ref{A}, this subgraph is a tree isomorphic to the graph of the
corresponding irreducible $U_{q_i}\widehat{\mathfrak
{sl}}_2$--module. Moreover, its root corresponds to an $i$--dominant
monomial. Therefore if a vertex of $\Omega_V(D)$ has no incoming
arrows of color $i$, then it corresponds to an $i$--dominant monomial.
In particular, if $m$ has no incoming arrows in $\Omega_V(D)$, then
$m$ is dominant. Since by our assumption $\chi_q(V)$ does not contain
any dominant monomials except for $m_+$, the graph $\Omega_V(D)$ is a
tree with root $m_+$.

Choose a sequence of weights $\la_1 > \la_2 > \ldots$ as above. We
prove by induction on $r$ the following statement $S_r$:

\medskip

The algorithm does not fail during the first $r$ steps. Let $\chi_r$
be the resulting polynomial after these steps. Then the coefficient
of each monomial $m$ in $\chi_r$ is not greater than that in
$\chi_q(V)$ and the coefficients of monomials of weights
$\la_1,\dots,\la_r$ in $\chi_r$ and $\chi_q(V)$ are equal. 
Furthermore, there exists a decomposition $D$ of $\chi_q(V)$, such that 
monomials in $\chi_r$ can be identified with
vertices in $\Omega_V(D)$ in such a way that all outgoing arrows from
vertices with $\uqgg$--weights $\la_1,\dots,\la_r$ go to
vertices of $\chi_r$. Finally, the $j$-th coloring of a monomial
$m$ in $\chi_r$ is just the number of vertices of type $m$ in $\chi_r$
which have incoming arrows of color $j$ in $\Omega_V(D)$.

\medskip

The statement $S_0$ is obviously true. Assume that the statement $S_r$
is true for some $r\geq 0$. Recall that at the $(r+1)$st step we
expand all monomials of $\chi_r$ of weight $\la_{r+1}$.

Let $m$ be a monomial of weight $\la_{r+1}$ in $\chi_r$, which enters
with coefficient $s$ and coloring $\{s_i\}_{i\in I}$. 

Then the monomial $m$ enters $\chi_q(V)$ with coefficient $s$ as
well. Indeed, $\Omega_V(D)$ is a tree, so all vertices $m$ have
incoming arrows from vertices of larger weight. By the statement $S_r$
this arrows go to vertices corresponding to monomials in $\chi_r$.

Suppose that $s_j<s$ for some $j\in I$. Then $m$ is
$j$--dominant. Indeed, otherwise each vertex of type 
$m$ in $\Omega_V(D)$ has an incoming arrow of color $j$ coming from a
vertex of higher weight. Then by the last part of the statement $S_r$, $s_j=s$.

Therefore the monomial $m$ is admissible, and the algorithm does not
fail at $m$. 

Consider the expansion $j_m(\chi_r)$. Let $\mu$ be as in \Ref{mu
form}.  In the $j$-th decomposition of $\chi_q(V)$, $m$ corresponds to
a root of a tree whose vertices can be identified with monomials in
$\mu$. We fix such an identification. Then monomials in $\mu$ get
identified with vertices in $\Omega_V(D)$.

Let $v$ be the vertex in $\Omega_V(D)$, corresponding to a monomial $n$ 
in $\mu$. Denote the coefficient of $n$ in $\chi_r$ by $p$
and the coloring by $\{p_i\}_{i\in I}$. We have two cases:

a) $p_j=p$. Then the last part of the statement $S_r$ implies that the
vertex $v$ does not belong to $\chi_r$. We add the monomial $n$ to
$\chi_r$ and increase $p_j$ by one (we have already identified it with
$v$).

b) $p_j<p$. Then by $S_r$ there exists a vertex $w$ in $\chi_r$ of
type $n$ with no incoming arrows of color $j$. We change the
decomposition $D_j$ by switching the vertices $v$ and $w$ and identify
$n$ with the new $v$. We also increase $p_j$ by one. (Thus, in this
case we do not add $n$ to $\chi_r$.)

In both cases, the statement $S_{r+1}$ follows.

Since the set of weights of monomials occuring in $\chi_q(V)$ is
contained in a finite set $\la_1,\la_2,\ldots,\la_N$, the statement
$S_N$ proves the first part of the theorem.

\corref{no dominant monomials} then implies the second part of the
theorem.
\end{proof}

We plan to use the above algorithm to compute explicitly the
$q$--characters of the fundamental representations of $\uqg$ and to
obtain their decompositions under $\uqgg$.

\begin{remark}
There is a similar algorithm for computing the ordinary characters of
finite-dimensional $\g$--modules (equivalently,
$\uqgg$--modules). That algorithm works for those representations
(called miniscule) whose characters do not contain dominant weights
other than the highest weight (for other representations the algorthim
does not work). However, there are very few miniscule representations
for a general simple Lie algebra $\g$. In contrast, in the case of
quantum affine algebras there are many representations whose
characters do not contain any dominant monomials except for the
highest weight monomials (for example, all fundamental
representations), and our algorithm may be applied to them.\qed
\end{remark}

\section{The fundamental representations}\label{fundamental
representations}

In this section we prove several theorems about the irreducibility 
of tensor products of fundamental representations.

\subsection{Reducible tensor products of fundamental representations
and poles of $R$-matrices} In this section we prove that the
reducibility of a tensor product of the fundamental representations is
always caused by a pole in the $R$-matrix.

We say that a monomial $m$ has {\em positive lattice support with base
$a$} if $m$ is a product $Y_{i,aq^n}^{\pm1}$ with $n\geq 0$.

\begin{lem}\label{positive support of fund}
All monomials in $\chi_q(V_{\omega_i}(a))$ have positive lattice
support with base $a$.
\end{lem}

\begin{proof}
For $\uqsl$, the statement follows from the explicit formula \Ref{sl2
fund} for $\chi_q(V_1(a))$. The $q$--character of any irreducible
representation $V$ of $\uqsl$ is a subsum of a product of the
$q$--characters of $V_1(b)$'s. Moreover, this subsum includes the
highest monomial. Hence if the highest weight monomial of $\chi_q(V)$
has positive lattice support with base $a$, then so do all monomials in
$\chi_q(V)$.

Now consider the case of general $\uqg$. Suppose there exists a
monomial in $\chi = \chi_q(V_{\omega_i}(a))$, which does not have
positive lattice support with base $a$. Let $m$ be a highest among
such monomials (with respect to the partial ordering by weights).

By \corref{no dominant monomials}, the monomial $m$ is not
dominant. In other words, if we rewrite $m$ as a product of
$Y_{i,b}^{\pm 1}$, we will have at least one generator in negative
power, say $Y^{-1}_{i_0,b_0}$.

Write $\tau_{i_0}(\chi)$ in the form \Ref{restriction improved}.  The
monomial $\tau_{i_0}(m)$ can not be among the monomials $\{m_pN_p\}$,
since $m$ contains $Y^{-1}_{i_0,b_0}$. Therefore $\tau_{i_0}(m) =
m_{p_0}N_{p_0}\ol M_{r_0,p_0}$ for some $\ol{M}_{r_0,p_0} \neq 1$,
which is a product of factors $\ol A^{-1}_{i,c}$. Let $m_1$ be a
monomial in $\chi$, such that $\tau_{i_0}(m_1) = m_{p_0}N_{p_0}$. Then
by Lemma \ref{A}, $m=m_1M_{r_0,p_0}$, where $M_{r_0,p_0}$ is obtained
from $\ol M_{r_0,p_0}$ by replacing all $\ol A^{-1}_{i,c}$ with
$A^{-1}_{i,c}$.

By construction, the weight of $m_1$ is higher than the weight of $m$,
so by our assumption, $m_1$ has positive lattice support with base
$a$. But then $m_{p_0}$ also has positive lattice support with base
$a$. Therefore all monomials in $m_{p_0}(1+\sum_r \ol{M}_{r,p})$ have
positive lattice support with base $a$. This implies that
$M_{r_0,p_0}$, and hence $m=m_1M_{r_0,p_0}$, has positive lattice
support with base $a$. This is a contradiction, so the lemma is
proved.
\end{proof}

\begin{remark}\label{strictly positive}
{}From the proof of Lemma \ref{positive support of fund} is clear that
the only monomial in $\chi_q(V_{\om_i}(a))$ which contains
$Y_{j,aq^n}^{\pm 1}$ with $n=0$ is the highest weight monomial
$Y_{i,a}$.\qed
\end{remark}

Let $V$ be a $\uqg$--module with the $q$--character $\chi_q(V)$. Define
the oriented graph $\Gamma_V$ as follows. The vertices of $\Gamma_V$
are monomials in $\chi_q(V)$ with multiplicities. Thus, there are
$\on{dim} V$ vertices. We denote the monomial corresponding to a
vertex $\al$ by $m_\al$. We draw an arrow from the vertex $\al$ to the
vertex $\beta$ if and only if $m_\beta=m_\al A^{-1}_{i,x}$ for some
$i\in I$, $x\in\C^\times$.

If $V$ is an irreducible $\uqsl$--module, then the graph $\Gamma_V$ is
connected. Indeed, every irreducible $\uqsl$--module is isomorphic to
a tensor product of evaluation modules. The graph associated to each
evaluation module is connected according to the explicit formulas for
the corresponding $q$--characters (see formula (4.3) in
\cite{FR:char}). Clearly, a tensor product of two modules with
connected graphs also has a connected graph.

\begin{lem} 
Let $\al\in \Gamma_V$ be a vertex with no incoming arrows. Then
$m_\al$ is a dominant monomial.
\end{lem}
\begin{proof}
Let $\al$ contain $Y^{-1}_{i,b}$ for some $i\in I$,
$b\in\C^\times$. We write the restricted $q$--character
$\tau_{i}(\chi_q(V))$ in the form \Ref{restriction improved}, where
each $m_p(1+\sum_r \ol M_{r,p})$ is a $q$--character of an
irreducible $U_{q_{i}}\widehat{\sw}_2$ module.

The monomial $\tau_{i}(m)$ contains $Y^{-1}_{i,b}$ and therefore can
not be among the monomials $\{ m_pN_p \}$. But the graphs of
irreducible $\uqsl$--modules are connected. So we obtain that
$\tau_{i}(m)=\tau_{i}(A_{i,c}^{-1})\tau_{i}(m')$ for some monomial
$m'$ in $\chi_q(V)$, and some $c\in\C^\times$. By Lemma \ref{A}, we
have $m=A_{i,c}^{-1}m'$ which is a contradiction.
\end{proof}

Now \corref{no dominant monomials} implies:

\begin{cor}
The graphs of all fundamental representations are connected.
\end{cor}

Let a monomial $m$ have lattice support with base $a$. We call $m$
{\em right negative} if the factors $Y_{i,aq^k}$ appearing in $m$, for
which $k$ is maximal, have negative powers.

\begin{lem}\label{right negative}
All monomials in the $q$--character of the fundamental representation
\linebreak $V_{\om_i}(a)$, except for the highest weight monomial, are
right negative.
\end{lem}

\begin{proof}
Let us show first that from the highest weight monomial $m_+$ there is
only one outgoing arrow to the monomial $m_1=m_+
A^{-1}_{i,aq_i}$. Indeed, the weight of a monomial that is connected
to $m_+$ by an arrow has to be equal to $\omega_i - \al_j$ for some $j
\in I$. The restriction of $V_{\om_i}(a)$ to $\uqg$ is isomorphic to
the direct some of its $i$th fundamental representation $V_{\om_i}$
and possibly some other irreducible representations with dominant
weights less than $\omega_i$. However, the weight $\omega_i - \al_j$
is not dominant for any $i$ and $j$. Therefore this weight has to
belong to the set of weights of $V_{\om_i}$, and the multiplicity of
this weight in $V_{\om_i}(a)$ has to be the same as that in
$V_{\om_i}$. It is clear that the only weight of the form $\omega_i -
\al_j$ that occurs in $V_{\om_i}$ is $\om_i-\al_i$, and it has
multiplicity one. By \thmref{1}, this monomial must have the form
$m_1=m_+ A^{-1}_{i,aq_i}$.

Now, the graph $\Gamma_{V_{\om_i}(a)}$ is connected. Therefore each
monomial $m$ in $\chi_q(V_{\om_i}(a))$ is a product of $m_1$ and
factors $A_{j,b}^{-1}$. Note that $m_1$ is right negative and all
$A_{j,b}^{-1}$ are right negative (this follows from the explicit
formula \eqref{express}). The product of two right negative monomials
is right negative. This implies the lemma.
\end{proof}

\begin{remark}    \label{2ri}
It follows from the proof of the lemma that the rightmost factor of
each non-highest weight monomial occurring in
$\chi_q(V_{\omega_i}(a))$ equals $Y_{j,aq^n}^{-1}$, where $n\geq
2r_i$.\qed
\end{remark}

Recall the definition of the normalized $R$--matrix $\ol{R}_{V,W}(z)$
from \secref{normR}. The following theorem was conjectured, e.g., in
\cite{AK}.

\begin{thm}    \label{poles}
Let $\{ V_k \}_{k=1,\dots,n}$, where $V_k =V_{\om_{s(k)}}(a_k)$, be a
set of fundamental representations of $\uqg$. The tensor product $V_1
\otimes \ldots \otimes V_n$ is reducible if and only if for some
$i,j\in \{ 1,\ldots,n\}, i\neq j$, the normalized $R$--matrix
$\ol{R}_{V_i,V_j}(z)$ has a pole at $z=a_j/a_i$.
\end{thm}

\begin{proof}
The ``if'' part of the Theorem is obvious. Let us explain the case
when $n=2$.  Let $\sigma: V_1 \otimes V_2 \arr V_2 \otimes V_1$ be the
transposition. By definition of $\ol{R}_{V_1,V_2}(z)$, the linear map
$\sigma \circ \ol{R}_{V_1,V_2}(z)$ is a homomorphism of
$\uqg$--modules $V_1 \otimes V_2 \arr V_2 \otimes V_1$. Therefore if
$\ol{R}_{V_1,V_2}(z)$ has a pole at $z=a_2/a_1$, then $V_1 \otimes
V_2$ is reducible. It is easy to generalize this argument to general
$n$.

Now we prove the ``only if'' part.

If the product $V_1\otimes\dots\otimes V_n$ is reducible, then the
product of the $q$--characters $\prod_{i=1}^n\chi_q(V_i)$ contains a
dominant monomial $m$ that is different from the product of the
highest weight monomials. Therefore $m$ is not right negative and $m$
is a product of some monomials $m_i'$ from $\chi_q(V_i)$. Hence at
least one of the factors $m_i'=m_i$ must be the highest weight
monomial and it has to cancel with the rightmost $Y_{i,b}^{-1}$
appearing in, say, $m_j'$.

According to \lemref{positive support of fund}, $m_j'=m_j M$ where $M$
is a product of $A_{s,a_jq^n}^{-1}$. By our assumption, the maximal
$n_0$ occurring among $n$ is such that $a_jq^{n_0}=a_iq_i^{-1}$. Using
Lemma \ref{diagonal lemma} we obtain that one of the diagonal entries
of $\ol{R}_{V_i,V_j}$ has a factor $1/(1-a_ia_j^{-1}z)$, which can not
be cancelled. Therefore $\ol{R}_{V_i,V_j}$ has a pole at
$z=a_j/a_i$. This proves the ``only if'' part. Moreover, we see that
the pole necessarily occurs in a diagonal entry.
\end{proof}

\subsection{The lowest weight monomial}
Our next goal is to describe (see \propref{irreducible tensor product}
below) the possible values of the spectral parameters of the fundamental
representations for which the tensor product is reducible.

First we develop an analogue of the formalism of \secref{str} from the
point of view of the lowest weight monomials. Recall the involution $I
\arr I, i \arr \bar{i}$ from \secref{qaal}. According to
\thmref{ChP2}(3), there is a unique lowest weight monomial $m_-$ in
$\chi_q(V_{\om_i}(a))$, and its weight is $-\om_{\bar i}$.

\begin{lem}
The lowest weight monomial of $\chi_q(V_{\om_i}(a))$ equals $Y_{\bar
i,aq^{r^\vee h^\vee}}^{-1}$.
\end{lem}

\begin{proof}
By Lemma \ref{dom}, $m_-$ must be antidominant. Thus, by Lemma
\ref{positive support of fund}, $m_-=Y_{\bar i,aq^{n_i}}^{-1}$ for
some $n_i>0$.

Recall the automorphism $w_0$ defined in \Ref{w0}.  The module
$V_{\om_{\bar i}}(a)$ is obtained from $V_{\om_i}(a)$ by pull-back
with respect to $w_0$. From the interpretation of the $q$--character
in terms of the eigenvalues of $\Phi_i^\pm(u)$, it is clear that the
$q$--character of $V_{\om_{\bar i}}(a)$ is obtained from the
$q$--character of $V_{\om_i}(a)$ by replacing each $Y_{j,b}^{\pm 1}$
by $Y_{\bar j,b}^{\pm 1}$. Therefore we obtain: $n_i=n_{\bar i}$.

Consider the dual module $V_{\om_i}(a)^*$. By \thmref{ChP2}(3), its
highest weight equals $\om_{\bar i}$. Hence $V_{\om_i}(a)^*$ is
isomorphic to $V_{\om_{\bar i }}(b)$ for some $b \in \C^\times$. Since
$\uqg$ is a Hopf algebra, the module $V_{\om_i}(a)^\otimes
V_{\om_i}(a)^*$ contains a one--dimensional trivial
submodule. Therefore the product of the corresponding $q$--characters
contains the monomial $m=1$. According to Lemma \ref{right negative},
it can be obtained only as a product of the highest weight monomial in
one $q$--character and the lowest monomial in another. Therefore,
$b=aq^{\pm n_i}$.

In the same way we obtain that $V_{\om_{\bar i}}(a)^*$ is isomorphic
to $V_{\om_i}(aq^{\pm n_i})$.

{}From formula \Ref{ss} for the square of the antipode, we obtain that
the double dual, $V_{\om_i}(a)^{**}$, is isomorphic to
$V_{\om_i}(aq^{-2r^\vee h^\vee})$. Since $n_i>0$, we obtain that
$n_i=r^\vee h^\vee$.
\end{proof}

Having found the lowest weight monomial in the $q$--characters of the
fundamental representations, we obtain using \thmref{ChP2} the lowest
weight monomial in the $q$--character of any irreducible module.

\begin{cor}\label{lowest in arbitrary}
Let $V$ be an irreducible $\uqg$--module. Let the highest weight
monomial in $\chi_q(V)$ be \be m_+=\prod_{i\in
I}\prod_{k=1}^{s_k}Y_{i,a_k^{(i)}}.  \ee Then the lowest weight
monomial in $\chi_q(V)$ is given by \be m_-=\prod_{i\in
I}\prod_{k=1}^{s_k}Y^{-1}_{\bar i,a_k^{(i)}q^{r^\vee h^\vee}}.  \ee
\end{cor}

We also obtain a new proof of the following corollary, which has been
previously proved in \cite{CP}, Proposition 5.1(b):

\begin{cor}
\be
V_{\om_i}(a)^* \simeq V_{\om_{\bar i}}(aq^{-r^\vee h^\vee}).
\ee
\end{cor}

Now we are in position to develop the theory of $q$--characters based
on the lowest weight and antidominant monomials as opposed to the
highest weight and dominant ones.

\begin{prop}\label{anti 1}
The $q$--character of an irreducible finite-dimensional $\uqg$ module
$V$ has the form \be \chi_q(V)=m_-(1+\sum N_p), \ee 
where $m_-$ is the
lowest weight monomial and each $N_p$ is a monomial in $A_{i,c}$, $i\in I$,
$c\in\C^\times$ (i.e., it does not contain any factors
$A_{i,c}^{-1}$).
\end{prop}

\begin{proof}
First we prove the following analogue of formula
\eqref{form of char}:
$$
\chi_q(V) = m_-( 1 +
\sum_p N_p),
$$
where each $N_p$ is a monomial in $A_{i,c}^{\pm 1}$,
$c\in\C^\times$. The proof of this formula is exactly the same as the
proof of Proposition 3 in \cite{FR:char}. The rest of the proof is
completely parallel to the proof of Theorem \ref{1}.
\end{proof}

\begin{lem}\label{no antidominant monomials}
The only antidominant monomial of $q$--character of a fundamental
representation is the lowest weight monomial.
\end{lem}

\begin{proof}
The proof is completely parallel to the proof of Lemma \ref{no
dominant monomials}.
\end{proof}

\begin{lem}\label{antipositive support of fund}
All monomials in a $q$--character of a fundamental representation are
products $Y_{i,aq^n}^{\pm1}$ with $n \leq r^\vee h^\vee$.
\end{lem}
\begin{proof}
The proof is completely parallel to the proof of Lemma \ref{positive
support of fund}.
\end{proof}

The combination of Lemmas \ref{positive support of fund} and
\ref{antipositive support of fund} yields the following result.

\begin{cor}
Let the highest weight monomial $m_+$ of the $q$--character of an
irreducible $\uqg$--module $V$ be a product of monomials $m_+^{(i)}$
which have positive lattice support with bases $a_i$. Let $s_i$ be the
maximal integer $s$, such that $Y_{k,a_iq^s}$ is present in
$m_+^{(i)}$ for some $k\in I$. Then any monomial $m$ in $\chi_q(V)$
can be written as a product of monomials $m^{(i)}$, where each
$m^{(i)}$ is a product of $Y_{j,a_iq^n}$ with $n \in \Z, 0\leq n\leq
s_i + r^\vee h^\vee$
\end{cor}

\subsection{Restrictions on the values of spectral parameters of
reducible tensor products of fundamental representations}

It was proved in \cite{KS} that $V_{\omega_i}(a) \otimes
V_{\omega_j}(b)$ is irreducible if $a/b$ does not belong to a
countable set. As M. Kashiwara explained to us, one can show that this
set is then necessarily finite. The following proposition, which was
conjectured, e.g., in \cite{AK}, gives a more precise description of
this set.

\begin{prop}\label{irreducible tensor product}
Let $a_i \in \C$, $i=1,\dots,n$, and suppose that the tensor product
of fundamental representations $V_{\omega_{i_1}}(a_1) \otimes \ldots
\otimes V_{\omega_{i_n}}(a_n)$ is reducible. Then there exist $m \neq
j$ such that $a_m/a_j=q^k$, where $k \in \Z$ and $2\leq k\leq r^\vee
h^\vee$.
\end{prop}

\begin{proof}
If $V_{\omega_{i_1}}(a_1) \otimes \ldots \otimes
V_{\omega_{i_n}}(a_n)$ is reducible, then
$\chi_q(V_{\omega_{i_1}}(a_1)) \ldots \chi_q(V_{\omega_{i_n}}(a_n))$
should contain a dominant term other than the product of the highest
weight terms. But for that to happen, for some $m$ and $j$, there have
to be cancellations between some $Y_{p,a_mq^{n}}^{-1}$ appearing in
$\chi_q(V_{\omega_{i_m}}(a_m))$ and some $Y_{p,a_jq^{l}}$ appearing in
$\chi_q(V_{\omega_{i_j}}(a_j))$. These cancellations may only occur if
$a_m/a_j=q^{\pm k}, k \in \Z$, and $0 \leq k \leq r^\vee h^\vee$, by
Lemmas \ref{positive support of fund} and \ref{antipositive support of
fund}. Moreover, $k\geq 2$ according to Remark \ref{2ri}.
\end{proof}

Note that combining \thmref{poles} and \propref{irreducible tensor
product} we obtain

\begin{cor}
The set of poles of the normalized $R$--matrix
$\ol{R}_{V_{\om_i}(a),V_{\om_j}(a)}(z)$ is a subset of the set $\{ q^k
| k \in \Z, 2 \leq |k| \leq r^\vee h^\vee \}$.
\end{cor}

\subsection{The $q$--characters of the dual representations}
In this subsection we show a simple way to obtain the $q$--character
of the dual representation.

Recall that ${\mc K}$ is given by \Ref{Q}.
\begin{lem}\label{unique}
Let $\chi_1,\chi_2 \in {\mc K}$. Assume that all dominant
monomials in $\chi_1$ are the same as in $\chi_2$ (counted with
multiplicities). Then $\chi_1=\chi_2$.  
\end{lem}
\begin{proof}
Consider $\chi=\chi_1-\chi_2$. We have $\chi \in{\mc K}$ and
$\chi$ has no dominant monomials. Then $\chi=0$ by Lemma \ref{dom}.
\end{proof}

Note that the similar statement is true for antidominant monomials.

\begin{prop}\label{dual}
Let $V_{\om_i}(a)$ be a fundamental representation. Then the
$q$--character of the dual representation $V_{\om_i}(a)^* \simeq
V_{\om_{\bar i}}(aq^{-r^\vee h^\vee})$ is obtained from the
$q$--character of $V_{\om_i}(a)$ by replacing each $Y_{i,aq^n}^{\pm
1}$ by $Y_{i,aq^{-n}}^{\mp 1}$.
\end{prop}
\begin{proof}
Let $\chi_1=\chi_q(V_{\om_{\bar i}}(aq^{-r^\vee h^\vee}))$ and
$\chi_2$ is obtained from $\chi(V_{\om_i}(a))$ by replacing
$Y_{i,aq^n}^{\pm 1}$ by $Y_{i,aq^{-n}}^{\mp 1}$. Then $\chi_1$ and
$\chi_2$ are elements in ${\mc K}$ with the only dominant monomial
$Y_{\bar i, aq^{-r^\vee h^\vee}}$ by Corollary \ref{no dominant
monomials} and Lemma \ref{no antidominant monomials}. Therefore
$\chi_1=\chi_2$ by Lemma \ref{unique}.
\end{proof}

\begin{remark}\label{dominant-antidominant}
One can define a similar procedure for obtaining the $q$--character of
the dual to any irreducible $\uqg$--module $V$. Namely, by
\thmref{ChP2}, $\chi_q(V)$ is a subsum in the product of
$q$--characters of fundamental representations. In particular, any
monomial $m$ in $\chi_q(V)$ is a product of monomials $m^{(i)}$ from
the $q$--characters of these fundamental representations and
Proposition \ref{dual} tells us what to do with each $m^{(i)}$. This
procedure is consistent because $\chi_q((V \otimes W)^*) = \chi_q(V^*)
\cdot \chi_q(W^*)$.

Note that under this procedure the dominant monomials go to the
antidominant monomials and vice versa.  $\qed$
\end{remark}

\end{document}